\documentclass{article}
\usepackage{amsthm}
\usepackage{amssymb}
\usepackage{amsmath}
\usepackage{epsfig}
\usepackage[dvips]{color}
\newcommand{\RR}{{\cal R}}
\newcommand{\R}{\mathrm{R}}
\newcommand{\QQ}{{\cal Q}}

\newcommand{\dom}{\mbox{dom}}
\newcommand{\el}{\mbox{el}}

\newtheorem{theorem}{Theorem}[section]

\newtheorem{lemma}[theorem]{Lemma}
\newtheorem{conjecture}[theorem]{Conjecture}
\newcount\claimno
\def\mytextindent#1{\indent\indent\llap{\rm#1\enspace}\ignorespaces}
\def\myitem{\par\hangindent30pt\mytextindent}
\def\myclaim#1#2{
   \bigskip\noindent\rlap{\rm(#1)}\ignorespaces
   \hangindent=33pt\hskip30pt{\sl#2}\bigskip}
\def\newclaim#1#2{
   \global\advance\claimno by 1\relax
   \bigskip\noindent\rlap{\rm(\the\claimno)}\ignorespaces
   \global\expandafter\edef\csname CLAIMLABEL#1\endcsname{\the\claimno}\relax
   \hangindent=33pt\hskip30pt{\sl#2}\bigskip}
\def\refclaim#1{\csname CLAIMLABEL#1\endcsname}
\def\mylabel#1{{\label{#1}}}
\def\junk#1{}
\def\rt#1{{\color{red}#1}}
\def\rt#1{#1}
\let\ppar=\par

\begin{document}
\title{Three-coloring triangle-free graphs on surfaces II.  
$4$-critical graphs in a disk\thanks{Supported by grant GACR~201/09/0197 of Czech Science Foundation.}}
\author{%
     Zden\v{e}k Dvo\v{r}\'ak\thanks{Computer Science Institute (CSI) of Charles University,
           Malostransk{\'e} n{\'a}m{\v e}st{\'\i} 25, 118 00 Prague, 
           Czech Republic. E-mail: {\tt rakdver@iuuk.mff.cuni.cz}.
	   Supported by the Center of Excellence -- Inst. for Theor. Comp. Sci., Prague, project P202/12/G061 of Czech Science Foundation.}
 \and
     Daniel Kr{\'a}l'\thanks{Warwick Mathematics Institute, DIMAP and Department of Computer Science, University of Warwick, Coventry CV4 7AL, United Kingdom. E-mail: {\tt D.Kral@warwick.ac.uk}.}
 \and
        Robin Thomas\thanks{School of Mathematics, 
        Georgia Institute of Technology, Atlanta, GA 30332. 
        E-mail: {\tt thomas@math.gatech.edu}.
        Partially supported by NSF Grants No.~DMS-0739366 and DMS-1202640.}
}
\date{June 29, 2017}
\maketitle
\begin{abstract}
Let $G$ be a plane graph of girth at least five.  We show that if there exists
a $3$-coloring $\phi$ of a cycle $C$ of $G$ that does not extend
to a $3$-coloring of $G$, then $G$ has a subgraph $H$ on $O(|C|)$ vertices
that also has no $3$-coloring extending $\phi$.  
This is asymptotically best possible and improves
a previous bound of Thomassen.  
In the next paper of the series we will use this result
and the attendant theory to prove a generalization to graphs on surfaces
with several precolored cycles.
\end{abstract}

\section{Introduction}

This paper is a part of a series aimed at studying the $3$-colorability
of graphs on a fixed surface that are either triangle-free, or have their
triangles restricted in some way.
Historically the first result in this direction is the following
classical theorem of Gr\"otzsch~\cite{Gro}.

\begin{theorem}
\label{grotzsch}
Every triangle-free planar graph is $3$-colorable.
\end{theorem}

Thomassen~\cite{thom-torus,Tho3list,ThoShortlist}
found three reasonably simple proofs of this statement.
Recently, two of us, in joint work with Kawarabayashi~\cite{DvoKawTho}
were able to design a linear-time algorithm to $3$-color triangle-free
planar graphs, and as a by-product found perhaps a yet simpler proof
of Theorem~\ref{grotzsch}.  Another significantly different proof was given
by Kostochka and Yancey~\cite{koyan}.

The statement of Theorem~\ref{grotzsch}
cannot be directly extended to any surface other than the sphere.
In fact, for every non-planar surface $\Sigma$ there are infinitely many
$4$-critical triangle-free graphs that can be drawn in $\Sigma$.
(A graph is \emph{$4$-critical} if it is not $3$-colorable, but every
proper subgraph is.)
For instance, the graphs obtained from an odd cycle of length five or more
by applying Mycielski's
construction \cite[Section~8.5]{BonMur} have that property.
Thus an algorithm for testing $3$-colorability of triangle-free graphs
on a fixed surface will have to involve more than just testing the
presence of finitely many obstructions.

The situation is different for graphs of girth at least five
by another deep theorem of Thomassen~\cite{thom-surf}, the following.

\begin{theorem}
\mylabel{thm:thomgirth5}
For every surface $\Sigma$ there are only finitely many $4$-critical
graphs of girth at least five that can be drawn in $\Sigma$.
\end{theorem}

Thus the $3$-colorability problem on a fixed surface
has a polynomial-time algorithm
for graphs of girth at least five, but the presence of cycles of
length four complicates matters.
Let us remark that there are no $4$-critical graphs of girth at least five
on the projective plane and the torus~\cite{thom-torus} and
on the Klein bottle~\cite{tw-klein}.

The only non-planar surface for which the $3$-colorability problem
for triangle-free graphs is fully characterized is the projective plane.
Building on earlier work of Youngs~\cite{Youngs}, Gimbel and
Thomassen~\cite{gimbel} obtained the following elegant characterization.
A graph drawn in a surface is a {\em quadrangulation} if every face
is bounded by a cycle of length four.

\begin{theorem}
\mylabel{thm:gimtho}
A triangle-free graph drawn in the projective plane is $3$-colorable if and only
if it has no subgraph isomorphic to a non-bipartite
quadrangulation of the projective plane.
\end{theorem}

For other surfaces there does not seem to be a similarly nice characterization.
Gimbel and Thomassen~\cite[Problem~3]{gimbel} asked whether there is
a polynomial-time algorithm to test the $3$-colorability of triangle-free
graphs embeddable in a fixed surface.
In a later paper of this series we will resolve this question in the
affirmative.
The algorithm naturally breaks into two steps.
The first is when the graph is
a quadrangulation, except perhaps for a bounded number of larger faces
of bounded size, which will be allowed to be precolored.
In this case there is a simple topological obstruction to the existence
of a coloring extension based on the so-called ``winding number" of
the precoloring.
Conversely, if the obstruction is not present and the graph is highly
``locally planar", then we can show that the precoloring can be
extended to a $3$-coloring of the entire graph.
This can be exploited to design a polynomial-time algorithm.
With additional effort the algorithm can be made to run in linear time.

The second step covers the remaining case, when the graph has either many faces
of size at least five, or one large face, and the same holds for every
subgraph.
In that case we show that the graph is $3$-colorable.
That is a consequence of the following theorem~\cite{proof-part4}, which will form the
cornerstone of the series of our papers.

\begin{theorem}
\mylabel{thm:corner}
There exists an absolute constant $K$ with the following property.
Let $G$ be a  graph drawn in a surface $\Sigma$ of Euler genus $\gamma$
with no separating cycles of length at most four,  and let $t$ be
the number of triangles in $G$.
If $G$ is $4$-critical, 
then $\sum|f|\le K(t+\gamma)$,
where the summation is over all faces $f$ of $G$ of length at least five.
\end{theorem}

If $G$ has girth at least five, then $t=0$ and every face has length
at least five.
Thus Theorem~\ref{thm:corner} implies Theorem~\ref{thm:thomgirth5},
and, in fact, improves the bound given by the proof of 
Theorem~\ref{thm:thomgirth5} in~\cite{thom-surf}.
The fact that our bound in Theorem~\ref{thm:corner} is linear in the
number of triangles is needed in our solution~\cite{proof-havel} of a 
problem of Havel~\cite{conj-havel}, as follows.

\begin{theorem}
\mylabel{havel}
There exists an absolute constant $d$ such that if $G$ is a planar
graph and every two distinct triangles in $G$ are at distance at least $d$,
then $G$ is $3$-colorable.
\end{theorem}

Our technique is a refinement of the standard method of reducible configurations.
We show that every sufficiently generic graph $G$ (i.e., a graph that is large enough
and cannot be decomposed to smaller pieces along cuts simplifying the problem)
embedded in a surface contains one of a fixed list of subgraphs.   Each such configuration
enables us to obtain a smaller $4$-critical graph $G'$ with the property that every $3$-coloring of $G'$
corresponds to a $3$-coloring of $G$.  Furthermore, we perform the reduction in such a way
that a properly defined weight of $G'$ is greater or equal to the weight of $G$.
A standard inductive argument then shows that the weight of every $4$-critical graph is bounded,
which also restricts its size.  Unfortunately, this brief exposition hides a large number of technical
details that need to be dealt with.

In this paper, we introduce this basic technique and apply it to prove the following special
case of Theorem~\ref{thm:corner}.

\begin{theorem}
\mylabel{thm:main}
Let $G$ be a graph of girth at least five drawn in the plane,  
let $C$ be a cycle in $G$, and let $\phi$ be a $3$-coloring of $C$ that does
not extend to a $3$-coloring of $G$.  
Then there exists a subgraph $H$ of $G$ containing $C$ such that
$|V(H)|\le 1715|V(C)|$ and $H$ has no $3$-coloring extending $\phi$.
\end{theorem}

After we obtained a proof of Theorem~\ref{thm:main}, but before we wrote
it down and made it public, the first author and 
Kawarabayashi~\cite{dvkaw} generalized Theorem~\ref{thm:main}
to list-coloring.
Their proof is about as long as ours, but has the added advantage that
it replaces $1715$ by a much smaller constant.
However, we are proceeding with publication of our paper, because we
need the theory it develops for the proof of Theorem~\ref{thm:corner}
for graphs of girth at least five, which will appear in the next
paper of our series.
It is natural to ask whether an analogue of Theorem~\ref{thm:corner}
restricted to graphs of girth at least five holds in the list-coloring setting.
An affirmative answer would be implied by the following conjecture,
see~\cite{hyper} for details.  Luke Postle (private communication) believes he has
a proof of Conjecture~\ref{conj-sh}, which however has not yet been written down.

\begin{conjecture}\label{conj-sh}
For every integer $k\ge 5$, there exists an integer $K$ with the following property.
Let $G$ be a planar graph of girth at least five, let $C_1,C_2$ be two
cycles in $G$ of lengths at most $k$, and for every $v\in V(G)$
let $L(v)$ be a set such that $|L(v)|=1$ if $v\in V(C_1\cup C_2)$ and $|L(v)|\ge3$ otherwise.
If there exists no proper coloring $\phi$ of $G$ such that $\phi(v)\in L(v)$
for every $v\in V(G)$,
then $G$ has a subgraph $H$ on at most $K$ vertices such that
$C_1$ and $C_2$ are subgraphs of $H$ and there exists no proper coloring
$\psi$ of $H$ such that $\phi(v)\in L(v)$ for every $v\in V(H)$.
\end{conjecture}

In order to avoid  duplication of work in the next paper of the series
we state many of the auxiliary
results in this paper in the more general setting of graphs on surfaces.  For this purpose,
we require some definitions introduced in the following section.  In Section~\ref{sec-reduc},
we describe more precisely what we mean by a reducible configuration, its appearance
in the considered graph and its reduction.  In Section~\ref{sec-col}, we show that the reductions
preserve $3$-colorings.  In Section~\ref{sec-disch}, we give the discharging argument used to
show the existence of a reducible configuration.  In Section~\ref{sec-reductions}, we argue that
the reductions preserve the assumptions of the theorem.  In Section~\ref{sec-winners}, we analyze
the change of the weights during the reduction.  In Section~\ref{sec-disk}, we combine
the results to prove Theorem~\ref{thm:main}.  Finally, in Section~\ref{sec-summary} we prove a technical
result summarizing the conclusions of this paper that will be used in the next paper~\cite{proof-part3} of this series.

\section{Definitions}
\label{sec:def}

All graphs in this paper are finite and simple, with no loops or parallel edges.

A \emph{surface}
is a compact connected $2$-manifold with (possibly null) boundary.  
Each component of the boundary
is homeomorphic to the circle, and we call it a \emph{cuff}.  For non-negative integers $a$, $b$ and $c$,
let $\Sigma(a,b,c)$ denote the surface obtained from the sphere by adding $a$ handles, $b$ crosscaps and
removing the interiors of $c$ pairwise disjoint closed discs.  
A standard result in topology shows that
every surface is homeomorphic to $\Sigma(a,b,c)$ for some choice of $a$, $b$ and $c$.
Note that $\Sigma(0,0,0)$ is a sphere, $\Sigma(0,0,1)$ is a closed disk, $\Sigma(0,0,2)$ is a cylinder,
$\Sigma(1,0,0)$ is a torus, $\Sigma(0,1,0)$ is a projective plane and $\Sigma(0,2,0)$ is a Klein bottle.
The \emph{Euler genus} $g(\Sigma)$ of the surface $\Sigma=\Sigma(a,b,c)$ is defined as $2a+b$.
For a cuff $C$ of $\Sigma$, let $\widehat{C}$ denote an open disk with boundary $C$ disjoint from $\Sigma$, and let $\Sigma+\widehat{C}$ be
the surface obtained by gluing $\Sigma$ and $\widehat{C}$ together, that is, by closing $C$ with a patch.
Let $\widehat{\Sigma}=\Sigma+\widehat{C_1}+\ldots+\widehat{C_c}$, where $C_1$, \ldots, $C_c$ are the cuffs of $\Sigma$,
be the surface without boundary obtained from $\Sigma$ by patching all the cuffs.

Consider a graph $G$ embedded in the surface $\Sigma$; when useful, we identify $G$ with the topological
space consisting of the points corresponding to the vertices of $G$ and the simple curves corresponding
to the edges of $G$.  We say that the embedding is \emph{normal} if every cuff of $\Sigma$ is equal to a cycle in $G$,
and we call such a cycle a \emph{ring}.
Throughout the paper, all graphs are embedded normally.
A \emph{face} $f$ of $G$ is a maximal arcwise-connected subset of $\Sigma-G$.
We write $F(G)$ for the set of faces of $G$.
The boundary of a face is equal to a union of closed walks of $G$, which we call the \emph{boundary walks} of $f$.  

Consider a ring $R$.  If $R$ is a triangle and at most one vertex of $R$ has degree greater than two in $G$, we say that $R$
is a \emph{vertex-like ring}.  A ring with only vertices of degree two is \emph{isolated}.  For a vertex-like ring $R$ that is
not isolated, the \emph{main} vertex of $R$ is its vertex of degree greater than two.
A vertex $v$ of $G$ is a \emph{ring vertex} if $v$ is belongs to a ring (i.e., $v$ is drawn in the boundary of $\Sigma$), 
and $v$ is \emph{internal} otherwise.  
A cycle $K$ in $G$ is \emph{separating} or \emph{separates the surface}
if $\widehat{\Sigma}-K$ has at least two components, 
and $K$ is \emph{non-separating} otherwise.
A cycle $K$ is \emph{contractible} if there exists a closed disk $\Delta\subseteq \Sigma$ with boundary equal to $K$.
A cycle $K$ \emph{surrounds the cuff $C$} if $K$ is not contractible in $\Sigma$, but it is contractible in $\Sigma+\widehat{C}$.
We say that $K$ \emph{surrounds a ring $R$} if $K$ surrounds the cuff incident with $R$.

Let $G$ be a graph embedded in a surface $\Sigma$, let the embedding be
normal, and let $\cal R$ be the set of rings of this embedding.
In those circumstances we say that $G$ is a \emph{graph in $\Sigma$
with rings $\cal R$.}
Furthermore, some vertex-like rings are designated as \emph{weak vertex-like rings}.  
At this point, let us remark that weak vertex-like rings are a technical
device designed to deal with cutvertices in Theorem~\ref{thm:corner}.
They will not play any role in this paper, but we need to introduce them
in order to be able to formulate the lemmas in this paper in such a
way that they can be applied in the proof of Theorem~\ref{thm:corner}.

For a vertex-like ring $R$, we define the \emph{length} of $R$ as $|R|=0$ if $R$ is weak and $|R|=1$ otherwise.
For a ring $R$ that is not vertex-like, the \emph{length} $|R|$ of $R$ is the number of vertices of $R$.
For a face $f$, by $|f|$ we mean the sum of the lengths of the boundary walks of $f$ (in particular, if an edge
appears twice in the boundary walks, it contributes $2$ to $|f|$).

Let $G$ be a graph with rings $\cal R$.  Let $H=\bigcup {\cal R}$ and let $H'$ be a (not necessarily induced) subgraph of $G$ obtained from $H$
by, for each weak vertex-like ring $R$, removing the main vertex and one of the non-main vertices of $R$
(or by removing two vertices of $R$ if $R$ has no main vertex),
so that $H'$ intersects $R$ in exactly one non-main vertex.  A {\em precoloring} $\psi$ of $\cal R$ is a 
$3$-coloring of the graph $H'$.  
A precoloring of $\cal R$ {\em extends to a $3$-coloring of $G$}
if there exists a $3$-coloring $\phi$ of $G$ such that $\phi(v)=\psi(v)$ for every $v\in V(H')$.
The graph $G$ is {\em $\cal R$-critical} if $G\neq H$ and for every proper subgraph
$G'$ of $G$ that contains $H$, there exists a precoloring of ${\cal R}$ that extends
to a $3$-coloring of $G'$, but not to a $3$-coloring of $G$.  For a precoloring $\kappa$ of $\cal R$
the graph $G$ is {\em $\kappa$-critical} if $\kappa$ does not extend to a $3$-coloring of $G$,
but it extends to a $3$-coloring of every proper subgraph of $G$ that contains $\cal R$.
Let us remark that if $G$ is $\kappa$-critical for some $\kappa$, then it is $\cal R$-critical,
but the converse is not true (for example, consider a graph consisting of a single ring with two chords).
On the other hand, if $\kappa$ is a precoloring of the rings of $G$ that does not extend to a $3$-coloring of $G$, then
$G$ contains a (not necessarily unique) $\kappa$-critical subgraph.

\section{Reducible configurations}
\label{sec-reduc}

By a  plane graph we mean a graph $G$ drawn in the plane with no crossings.
Thus $G$ has exactly one unbounded face, called the {\em infinite face};
all the other faces are called {\em finite}.
An {\em isomorphism} of plane graphs maps finite faces to finite faces and
the infinite face to the infinite face.

A {\em configuration} is a quintuple $\gamma=(G,{\cal F},d,{\cal I}, {\cal A})$, where
{
\myitem{$\bullet$}$G$ is a plane graph,
\myitem{$\bullet$}$\cal F$ is a set of finite faces of $G$,
\myitem{$\bullet$}$d$ is a function that maps a set $\dom(d)\subseteq V(G)$
to $\{3,4,\ldots\}$,
\myitem{$\bullet$}$\cal I$ is a subset of $V(G)\setminus\dom(d)$, and
\myitem{$\bullet$}$\cal A$ is a subset of $V(G)\setminus\dom(d)$ of size zero or two.
}

\noindent
If $\gamma$ is a configuration, then we write $G_\gamma:=G$, 
${\cal F}_\gamma:={\cal F}$, $d_\gamma:=d$, ${\cal I}_\gamma:={\cal I}$
and ${\cal A}_\gamma:={\cal A}$.

\begin{figure}
\epsfbox{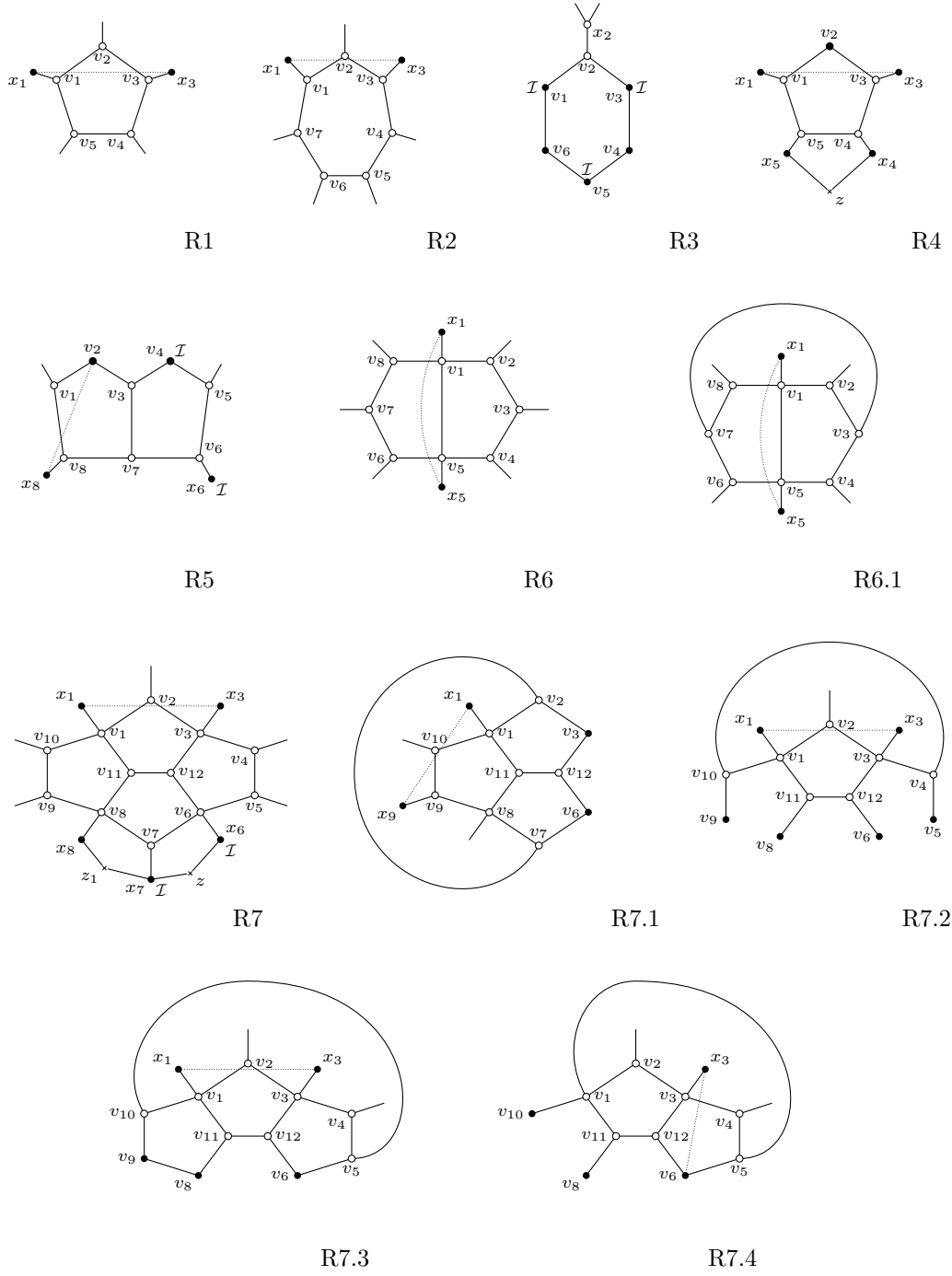}
\caption{Reducible configurations.}
\label{fig-redu}
\end{figure}

Two configurations $\gamma$ and $\gamma'$ are {\em isomorphic} if there exists
an isomorphism $\phi$ of the plane graphs $G_\gamma$ and $G_{\gamma'}$ that 
maps ${\cal F}_\gamma$ to ${\cal F}_{\gamma'}$,
${\cal I}_\gamma$ to ${\cal I}_{\gamma'}$,
${\cal A}_\gamma$ to ${\cal A}_{\gamma'}$,
$\dom(d_\gamma)$ to $\dom(d_{\gamma'})$ and
$d_\gamma(v)=d_{\gamma'}(\phi(v))$ for every $v\in \dom(d_\gamma)$.
Figure~\ref{fig-redu} contains the depictions of several configurations,
using the following conventions.
The graph $G_\gamma$ is drawn in the figure (ignoring the ``half-edges" 
and dashed edges for a moment);
 ${\cal F}_\gamma$ consists of all the finite faces of $G_\gamma$ that
do not include any half-edges in their interior;
the elements of  ${\cal I}_\gamma$ are indicated by ${\cal I}$
next to them;
if ${\cal A}_\gamma$ is non-empty, then the two vertices of ${\cal A}_\gamma$ are joined by a dashed edge;
the set $\dom(d_\gamma)$ consists of vertices drawn by empty circles;
and the value $d_\gamma(v)$ is equal to the
number of edges and half-edges incident with $v$ in the figure.
A configuration is {\em good} if it is isomorphic to one of the
configurations depicted in Figure~\ref{fig-redu}.

Let $\gamma$ be a good configuration and either let $H=G_\gamma$, or let
$H$ be a plane graph obtained from $G_\gamma$ by identifying
two vertices of $V(G_\gamma)\setminus\dom(d_\gamma)$ that are at distance
at least five in $G_\gamma$.
(The latter is only possible when $\gamma$ is $\R7$ or $\R7.2$.)
In those circumstances we say that $H$ is an {\em imprint} of $\gamma$.
It follows that every face in ${\cal F}_\gamma$ may be regarded as a face
of $H$, and that $\dom(d_\gamma)\subseteq V(H)$.

Let $G$ be a graph in a surface $\Sigma$ with rings $\cal R$.
We say that a configuration $\gamma$ {\em faintly appears} in \rt{$G$} if
{
\myitem{$\bullet$} some imprint $H$ of $\gamma$ is a subgraph of $G$,
\myitem{$\bullet$} every face in ${\cal F}_\gamma$ is a face of $G$,
\myitem{$\bullet$} $\dom(d_\gamma)\cap V({\cal R})=\emptyset$,
\myitem{$\bullet$} if $v\in\dom(d_\gamma)$, then $\deg_{G}(v)=d_\gamma(v)$, and
\myitem{$\bullet$} at most one vertex of ${\cal I}_\gamma$ belongs to $V({\cal R})$.
}

\noindent 
If a configuration $\gamma$ faintly appears in \rt{$G$}, then
we say that a subgraph $J$ of $G$ {\em touches} $\gamma$ if
an edge of $J$ is incident with a face in ${\cal F}_\gamma$.
We say that $\gamma$ {\em weakly appears} in \rt{$G$}
if it faintly appears and
{
\myitem{$\bullet$} no cycle of length at most four distinct from rings touches $\gamma$ and
if $\gamma$ is $\R7$, then $x_3\neq x_7$ or $x_1\neq x_6$,
\myitem{$\bullet$} if $u,v\in \dom(d_\gamma)$ are adjacent in $G$,
then $u,v$ are adjacent in $G_\gamma$, 
\myitem{$\bullet$} if $\gamma$ is isomorphic to $\R4$ and the
vertices corresponding to $x_4$ and $x_5$ both belong to $\cal R$,
then the vertex corresponding to $v_2$ does not belong to $\cal R$. 
}

Let a good configuration $\gamma$ weakly appear in $G$.
We wish to define a new graph $G'$ in $\Sigma$ with rings ${\cal R}$.
For the definition we need to distinguish several cases.
Assume first that $\gamma$ is not isomorphic to $\R4$.
Let the graph $G'$ be obtained
from $G\backslash \dom(d_\gamma)$ 
by adding an edge joining the vertices in 
${\cal A}_\gamma$ if ${\cal A}_\gamma\neq\emptyset$,
and by identifying the vertices in ${\cal I}_\gamma$.  If parallel edges are created,
remove all edges but one from each bunch of parallel edges, so that each edge of $G'$
corresponds to a unique edge of $G$.
Since no cycle of length at most four touches $\gamma$
and if $\gamma$ is $\R7$, then $x_3\neq x_7$ or $x_1\neq x_6$,
it follows that $G'$ has no loops.
It also follows that $\cal R$ is a set of rings for $G'$.
We will refer to the added edge as the {\em new edge} and to the vertex
that resulted from the identification of vertices as the {\em new vertex}.
If two vertices $u,v\in{\cal I}_\gamma$ have a common neighbor $x\in V(G_\gamma)\setminus \dom(d_\gamma)$
and $w$ is the new vertex arising by identification of $u$ and $v$, then we call
the edge $wx$ {\em squashed}.

We also need to specify an embedding of $G'$ in $\Sigma$. 
There is a unique natural way to make the edge additions and
vertex identifications inside the faces of ${\cal F}_\gamma$,
and that is how the embedding of $G'$ will be defined.
Formally, for every pair $u,v\in{\cal A}_\gamma$
and every pair $u,v\in{\cal I}_\gamma$
of distinct vertices
we define the {\em replacement $u,v$-path} as the shortest
path from $u$ to $v$ in $G_\gamma$.
It follows by inspecting all the good configurations
that the replacement path is unique.
Now we identify $u$ and $v$ or join them by an edge 
along the replacement $u,v$-path $P$,
with the proviso that if $P$ includes a vertex 
$v\in V(G_\gamma)\setminus\dom(d_\gamma)$
(specifically, vertex $v_4$ or $v_6$ of $\R3$ or vertex $z$ of $\R7$),
then prior to making the edge addition or vertex identification
we shift $P$ slightly into the unique face $f$ of ${\cal F}_\gamma$ incident with $v$.
Note that $P$ stays in $\Sigma$ and its homotopy does not change by such a shift.
This completes the definition of $G'$ when $\gamma$ is not $\R4$.

Now let $\gamma$ be $\R4$.
If not both $x_4$ and $x_5$ belong to $\cal R$, then we proceed as
above, treating the configuration as if $\{x_4,x_5\}$ belonged to 
${\cal I}_\gamma$; that is, identifying those vertices.
We may therefore assume that both $x_4,x_5$ belong to ${\cal R}$.
Let $\phi$ be a $3$-coloring of $\cal R$; the definition of $G'$ will
now depend on $\phi$.
If $\phi(x_4)=\phi(x_5)$, then we define $G'$ exactly as in the
previous two paragraphs; in particular, we do not identify $x_4$ and $x_5$.
If $\phi(x_4)\ne\phi(x_5)$, then we let $G'$ be obtained from
$G\backslash \{v_1,v_3,v_4,v_5\}$ by identifying $v_2$ and $x_5$
along the ``replacement path" $v_2v_1v_5x_5$ (we do not add the edge between
$x_1,x_3\in {\cal A}_\gamma$).
Let us remark that the last condition in the definition of 
weak appearance guarantees that in this case $v_2$ does not belong to $\cal R$.
Then $G'$ is a graph in $\Sigma$ with rings ${\cal R}$,
and we say that it is the {\em $\gamma$-reduction of $G$.}
When we wish to emphasize the dependence on $\phi$ we will say that
$G'$ is the {\em $\gamma$-reduction of $G$ with respect to $\phi$.}

\section{Colorings}
\label{sec-col}

In this section, we show that each $3$-coloring of the $\gamma$-reduction
of a  graph $G$ extends to a $3$-coloring of $G$.
Most of the reductions were used earlier~\cite{Gro,thom-torus},
but $\R5$, $\R7$ and their variants seem to be new.
For the sake of completeness we include proofs of extendability
for all good configurations.

\begin{lemma}\label{lemma-color}
Let $G$ be a graph in a surface $\Sigma$ with rings $\cal R$, 
let $\gamma$ be a good configuration that weakly appears in $G$,
let $\phi_0$ be a $3$-coloring of $\cal R$,
and let $G_1$ be the $\gamma$-reduction of $G$ with respect to $\phi_0$.
If $\phi_0$ extends to a $3$-coloring of $G_1$, then it extends
to a $3$-coloring of $G$.
\end{lemma}

\proof
Let $\gamma$ be as stated, and let the vertices of $G_\gamma$ be
labeled as in Figure~\ref{fig-redu}.
Let $\phi$ be a $3$-coloring of $G_1$ that extends the coloring $\phi_0$.
Then $\phi$ can be regarded as a $3$-coloring of 
$G\setminus\dom(d_\gamma)$, and our objective is to extend it
to a $3$-coloring of $G$.
For each vertex $v_i\in \dom(d_\gamma)$ that has a unique
neighbor outside of the configuration, let this neighbor be denoted by $x_i$.
We will use the following easy observations: 

{
\newclaim{cl-xext}{Suppose that $u_1,u_2\in V(G)$ are adjacent vertices 
   of degree three,
   $w_1$ and $w_2$ are the neighbors of $u_1$ distinct from $u_2$ and 
   $w_3$ and $w_4$ are the neighbors of $u_2$ distinct from $u_1$.
   A $3$-coloring $\psi$ of $w_1$, \ldots, $w_4$ extends to $u_1$ and 
   $u_2$, unless $\psi(w_1)=\psi(w_3)\neq\psi(w_2)=\psi(w_4)$
   or $\psi(w_1)=\psi(w_4)\neq\psi(w_2)=\psi(w_3)$.}
}

{
\newclaim{cl-pext}{Let $P=u_1u_2\ldots u_k$ be a path in $G$  and 
   $L_1$, \ldots, $L_k$ lists of colors of size two,
   such that $L_i\neq L_j$ for some $1\le i<j\le k$.  
   Then there exist colorings $\psi_1$, $\psi_2$ and $\psi_3$ of $P$ such 
   that $\psi_i(v_j)\in L_j$ for $1\le i\le 3$ and $1\le j\le k$, and
   for each $1\le i<j\le 3$ either    
   $\psi_i(u_1)\neq \psi_j(u_1)$ or $\psi_i(u_k)\neq \psi_j(u_k)$.}
}

Let us now consider each configuration separately.
\medskip

\noindent{\bf Configurations $\R1$ and $\R2$.} 
Each of the vertices of the cycle $v_1v_2\ldots v_k$ 
(where $k=5$ for the configuration $\R1$ and $k=7$ for $\R2$) 
has a list of two available colors, and the lists of $v_1$ and $v_3$ 
are not the same.
By (\refclaim{cl-pext}), there exists a coloring of the path 
$v_1\ldots v_k$ from these lists such that the colors
of $v_1$ and $v_k$ are not the same, giving a coloring of $G$, as desired.
\medskip

\noindent{\bf Configuration $\R3$.} 
The vertices $v_1$, $v_3$ and $v_5$ inherit the color of the new vertex.
Then we can color the vertices $x_2$ and $v_2$ in order,
because at the time each of those vertices is colored
it is adjacent to vertices of at most two different colors.
\medskip

\noindent{\bf Configuration $\R4$.}
Suppose first that at least one of $x_4$ and $x_5$ is internal, or 
that both belong to $\cal R$ and $\phi_0(x_4)=\phi_0(x_5)$.
If $\phi(x_1)=\phi(v_2)$, then color the vertices in the order 
$v_3$, $v_4$, $v_5$ and $v_1$ (each of them has neighbors of at most two 
different colors when it is being colored).  
The case that $\phi(x_3)=\phi(v_2)$ is symmetric.  
Therefore, we may assume that $\phi(x_1)=1$, $\phi(v_2)=2$
and $\phi(x_3)=3$.  
Set $\phi(v_1)=3$ and $\phi(v_3)=1$ and extend the coloring to 
$v_4$ and $v_5$ by (\refclaim{cl-xext}).
Then $\phi$ is a desired $3$-coloring of $G$.

We may therefore assume that both $x_4$ and $x_5$ belong to $\cal R$
and $\phi_0(x_4)\neq\phi_0(x_5)$.
In this case, the definition of $\gamma$-reduction 
ensures that $\phi(v_2)=\phi(x_5)$.  
We may assume that $\phi(v_2)=\phi(x_5)=1$ and $\phi(x_4)=2$.  
Let us set $\phi(v_4)=1$ and color $v_3$, $v_1$ and $v_5$ in this order.
\medskip

\noindent{\bf Configuration $\R5$.} 
The reduction ensures that $\phi(v_2)\neq \phi(x_8)$ and 
$\phi(v_4)=\phi(x_6)$.
If $\phi(v_2)=\phi(v_4)$, then $\phi$ extends---color the vertices 
in the order $v_1$, $v_8$, $v_5$, $v_6$, $v_7$ and $v_3$, and observe that
for each of these vertices, at most two different colors appear on 
already colored neighbors.  
Thus we may assume that $\phi(v_2)=1$ and $\phi(v_4)=\phi(x_6)=2$. 
We set $\phi(v_3)=3$ and  
$\phi(v_7)=2$, and color the vertices $v_5$ and $v_6$ by (\refclaim{cl-xext}).  
As $\phi(x_8)\neq\phi(v_2)\neq \phi(v_7)$, the observation (\refclaim{cl-xext}) 
implies that the coloring extends to $v_1$ and $v_8$.
\medskip

\noindent{\bf Configurations $\R6$ and $\R6.1$.} 
In both cases, the reduction ensures that $\phi(x_1)\neq \phi(x_5)$, 
say $\phi(x_1)=1$ and $\phi(x_5)=2$.  
If $\phi(x_6)=1$, then set $\phi(v_5)=1$, and color the vertices 
in order $v_4$, $v_3$, $v_2$, $v_1$, $v_8$, $v_7$ and $v_6$.  
Therefore, we may assume that this is not the case.
By symmetry, we may also assume that $\phi(x_4)\neq 1$ and 
$\phi(x_2),\phi(x_8)\neq 2$.
If $\phi(x_2)=\phi(x_8)=3$, then set $\phi(v_1)=3$, $\phi(v_5)=1$ and 
color $v_6$, $v_7$, $v_8$, $v_4$, $v_3$ and $v_2$ in this order.  
Otherwise, by symmetry we may assume that $\phi(x_2)=1$.
If $v_3$ and $v_7$ are adjacent, or if $\phi(x_3)\neq 1$, 
then set $\phi(v_3)=\phi(v_5)=1$ and color 
$v_4$, $v_6$, $v_7$, $v_8$, $v_1$ and $v_2$ in this order.  
Therefore, assume that $v_3$ and $v_7$ are not adjacent and $\phi(x_3)=1$.

If $\phi(x_6)=3$, then set $\phi(v_4)=1$, $\phi(v_1)=\phi(v_3)=2$ 
and $\phi(v_2)=\phi(v_5)=3$ and color $v_8$, $v_7$ and $v_6$ in this order.  
Thus, assume that $\phi(x_6)=2$.
By the argument symmetrical to the one used for $x_3$, we conclude that 
$\phi$ extends unless $\phi(x_7)=2$.
If $\phi(x_8)=3$, then set $\phi(v_4)=\phi(v_6)=\phi(v_8)=1$, 
$\phi(v_1)=\phi(v_3)=2$ and $\phi(v_2)=\phi(v_5)=\phi(v_7)=3$.  
Thus assume that $\phi(x_8)=1$ and by symmetry, $\phi(x_4)=2$.  
In this case, set $\phi(v_5)=\phi(v_7)=1$, $\phi(v_1)=\phi(v_3)=2$ and
$\phi(v_2)=\phi(v_4)=\phi(v_6)=\phi(v_8)=3$.
\medskip

\noindent{\bf Configuration $\R7$.}
The reduction ensures that $\phi(x_1)\neq \phi(x_3)$, say $\phi(x_1)=1$ and $\phi(x_3)=2$. To preserve the symmetry of the configuration, let us for a while ignore the identification of $x_6$ and $x_7$.

Suppose first that $\phi(x_8)=1$.  By (\refclaim{cl-pext}), there exists a coloring $\psi$ of the path $v_1v_2\ldots v_8$ such that $\psi(v_1)=\psi(v_8)\in\{2,3\}$.  We can extend $\psi$ to $v_{12}$ and $v_{11}$.  By (\refclaim{cl-xext}), if $\phi(x_9)\neq\phi(x_{10})$, then $\psi$ extends to $v_9$ and $v_{10}$ as well.
We next analyze the case that $\phi(x_9)=\phi(x_{10})=c$. Set $\phi(v_{11})=1$.  If $\phi(x_2)=1$, then color $v_3$ by $1$, and color the vertices $v_4$, $v_5$, \ldots, $v_{10}$, $v_1$, $v_2$, $v_{12}$ in this order. If $\phi(x_6)=2$, then color $v_{12}$ by $2$ and extend the coloring to the $10$-cycle $v_1\ldots v_{10}$.
Therefore, assume that $\phi(x_2)\neq 1$ and $d=\phi(x_6)\neq 2$.
Let us distinguish several cases:
\begin{itemize}
\item $d=3$, $\phi(x_4)=1$ and $\phi(x_5)=3$: In this case, set $\phi(v_{12})=3$, $\phi(v_3)=1$ and color $v_2$, $v_1$, $v_{10}$, $v_9$, \ldots, $v_4$ in order.
\item $d=1$ and $\phi(x_4)=\phi(x_5)$: Set $\phi(v_3)=1$ and color the vertices $v_2$, $v_1$, $v_{10}$, $v_9$, \ldots, $v_6$, $v_{12}$ in order.  Note that $\phi(v_3)=1\neq\phi(v_6)$, thus $\phi$ extends the coloring to $v_4$ and $v_5$ by (\refclaim{cl-xext}).
\item Otherwise, set $\phi(v_2)=1$, $\phi(v_3)=3$, $\phi(v_{12})=2$, $\phi(v_6)=4-d$, and color vertices $v_7$, \ldots, $v_{10}$, $v_1$ in order.  By (\refclaim{cl-xext}), this coloring extends to $v_4$ and $v_5$.
\end{itemize}
We conclude that if $\phi$ does not extend to the empty-circle vertices, then $\phi(x_8)=c_1\neq 1$,
and by the symmetry, $\phi(x_6)=c_2\neq 2$.

There are four possible colorings of $v_1$ and $v_8$ (two choices of colors for each of these vertices, so that the color of $v_1$ is not $1$ and the color of $v_8$ is not $c_1$). By (\refclaim{cl-xext}), out of these four colorings, all but at most one extend to $v_9$ and $v_{10}$; if such a coloring of $v_1$ and $v_8$ exists, let it be denoted by $\omega_1$; otherwise, set $\omega_1(v_1)=1$ and $\omega_1(v_8)=c_1$. Symmetrically, let $\omega_2$ be the unique coloring of $v_3$ and $v_6$ such that $\omega_2(v_3)\neq 2$, $\omega_2(v_6)\neq c_2$ and $\omega_2$ does not extend to $v_4$ and $v_5$, if such a coloring exists, and $\omega_2(v_3)=2$ and $\omega_2(v_6)=c_2$ otherwise.

If $\phi(x_2)=2$, then let $a=2$, otherwise let $a=3$.  Note that any color $c\neq 2=\phi(x_3)$ satisfies $|\{a,c,\phi(x_2)\}|=2$.  In the following cases, we can extend $\phi$ to a coloring $\psi$ of the path $v_1v_{10}v_9v_8v_7v_6$ such that $\psi(v_1)=a$ and $b=\psi(v_6)\neq\omega_2(v_6)$:
\begin{itemize}
\item $\omega_1(v_1)\neq a$: choose $b\not\in\{\phi(x_6), \omega_2(v_6)\}$, color $v_7$ and $v_8$, and note that we can extend this coloring to $v_9$ and $v_{10}$ by the definition of $\omega_1$.
\item $\omega_2(v_6)=c_2$: color the vertices $v_{10}$, $v_9$, \ldots, $v_6$ in this order.
\item $\phi(x_7)\not\in\{c_1,\omega_1(v_8)\}\cap\{c_2,\omega_2(v_6)\}$ or $\{c_1,\omega_1(v_8)\}=\{c_2,\omega_2(v_6)\}$:
excluding the previous two cases, we may assume that $c_1\neq \omega_1(v_8)$ and $c_2\neq\omega_2(v_6)$.  Color $v_8$ by the color $d\not\in\{c_1,\omega_1(v_8)\}$ and $v_6$ by the color $b\not\in\{c_2,\omega_2(v_6)\}$, extend the coloring to $v_9$ and $v_{10}$ by the definition of $\omega_1$, and observe that $|\{\phi(x_7),b,d\}|\le 2$,
thus $v_7$ can be colored as well.
\end{itemize}

If such a coloring $\psi$ exists, then choose a color $c\neq\phi(x_3)$ such that $c=b$ or $\{b,c\}\neq\{a, \psi(v_8)\}$; this ensures that the coloring extends to $v_{11}$ and $v_{12}$ by (\refclaim{cl-xext}). Since $b\neq \omega_2(v_6)$, this coloring extends to $v_4$ and $v_5$ as well.  Finally, the choice of $a$ ensures that $|\{a,c,\phi(x_2)\}|=2$, hence the coloring extends to $v_2$. Therefore, we may assume that such the coloring $\psi$ does not exist, i.e., $\omega_1(v_1)=a$, $\omega_2(v_6)\neq c_2$, $\{c_1,\omega_1(v_8)\}\neq \{c_2,\omega_2(v_6)\}$ and
$\phi(x_7)\in\{c_1,\omega_1(v_8)\}\cap\{c_2,\omega_2(v_6)\}$.

Let us now distinguish two cases:
\begin{itemize}\item $\phi(x_9)\neq \phi(x_{10})$:  By (\refclaim{cl-xext}), $a=\omega_1(v_1)=\phi(x_9)$. If $c_1\neq a$, then set $\phi(v_1)=\phi(v_8)=a$ and color $v_{10}$, $v_9$, $v_7$, $v_6$, \ldots, $v_2$ in this order ($v_2$ can be colored by the choice of $a$), and color $v_{12}$ and $v_{11}$;
hence, assume that $c_1=a$.

If $\phi(x_{10})=5-a$, then set $\phi(v_1)=\phi(v_8)=5-a$, $\phi(v_{10})=a$, and $\phi(v_9)=1$.  Note that $\phi(x_7)\in\{c_1,\omega_1(v_8)\}=\{a,5-a\}$ and $\{c_2,\omega_2(v_6)\}=\{1,\phi(x_7)\}$.  Set $\phi(v_7)=1$ and choose $\phi(v_6)\not\in\{c_2,\omega_2(v_6)\}$, i.e., $\phi(v_6)=5-\phi(x_7)$.  Extend the coloring to $v_2$, $v_3$, $v_{12}$ and $v_{11}$ in this order.  As $\phi(v_6)\neq\omega_2(v_6)$, this coloring extends to $v_4$ and $v_5$, giving a coloring of the whole configuration.

Therefore, assume that $\phi(x_{10})=1$.  Then $\omega_1(v_8)=1$ and $\phi(x_7)\in\{1, a\}$. Let us set $\phi(v_1)=\phi(v_7)=\phi(v_9)=5-a$, $\phi(v_{10})=a$ and $\phi(v_8)=1$. Let us choose color $\phi(v_6)\not\in\{c_2,\omega_2(v_6)\}$; note that $\phi(v_6)\neq 5-a$, since $\{c_2,\omega_2(v_6)\}\neq \{c_1,\omega_1(v_8)\}=\{1,a\}$.  Color $v_2$ and $v_3$, and extend the coloring to $v_4$ and $v_5$ (this is possible, since $\phi(v_6)\neq\omega_2(v_6)$).  We may assume that this coloring does not extend to $v_{11}$ and $v_{12}$, i.e., $\{\phi(v_3),\phi(v_6)\}=\{1,5-a\}$, hence $\phi(v_3)=5-a$ and $\phi(v_6)=1$.  As $\phi(v_6)\not\in\{c_2,\omega_2(v_6)\}$, we get $\{c_2,\omega_2(v_6)\}=\{a,5-a\}$ and $\phi(x_7)=a$.  Since $c_2\neq 2$, we have $c_2=3$ and$\omega_2(v_6)=2$.  As $\omega_2(v_3)\neq 2$, it follows that $\phi(x_4)=2$ and $\phi(x_5)\neq 2$.

Consider the coloring $\psi$ with $\psi(v_8)=5-a$, $\psi(v_7)=\psi(v_9)=1$, $\psi(v_6)=2$, $\psi(v_3)=\psi(v_5)=4-\phi(x_5)$ and $\psi(v_4)=\phi(v_5)$, and assume that this coloring does not extend to the coloring of the whole configuration.  On one hand, we may color $v_1$ by $a$ and $v_{10}$ by $5-a$; then $\psi$ extends to $v_2$ by the definition of $a$, and since it does not extend to $v_{11}$ and $v_{12}$, we have $\{a,5-a\}=\{2,4-\phi(x_5)\}$, and $\phi(x_5)=1$. On the other hand, we may color $v_1$ by $5-a$, $v_{12}$ by $1$ and $v_{10}$ and $v_{11}$ by $a$. Since this coloring does not extend to $v_2$, we have $|\{5-a,3,\phi(x_2)\}|=3$, and $a=3$ and $\phi(x_2)=1$.  In that case, we can color the configuration by setting $\phi(v_3)=\phi(v_6)=\phi(v_8)=1$, $\phi(v_1)=\phi(v_5)=\phi(v_7)=\phi(v_9)=\phi(v_{12})=2$ and $\phi(v_2)=\phi(v_4)=\phi(v_{10})=\phi(v_{11})=3$.
\item $\phi(x_9)=\phi(x_{10})$:  By symmetry, we may also assume that $\phi(x_4)=\phi(x_5)$. At this point, we use the second relation guaranteed by the reduction, $\phi(x_7)=c_2$. If $c_2\neq 3$, then set $\phi(v_7)=3$, $\phi(v_8)=1$ and $\phi(v_6)=2$, color the $5$-cycle $v_1v_2v_3v_{12}v_{11}$, and extend the coloring to $v_4$, $v_5$, $v_9$ and $v_{10}$ by (\refclaim{cl-xext}).  Thus, we can assume that $c_2=3$.

If $\phi(x_2)\neq 1$, then set $\phi(v_2)=\phi(v_6)=\phi(v_8)=1$, $\phi(v_1)=\phi(v_7)=\phi(v_{12})=2$ and $\phi(v_3)=\phi(v_{11})=3$, and extend the coloring to $v_4$, $v_5$, $v_9$ and $v_{10}$ by (\refclaim{cl-xext}).

Finally, if $\phi(x_2)=1$, then set $\phi(v_2)=\phi(v_8)=5-c_1$, $\phi(v_1)=c_1$, $\phi(v_3)=\phi(v_7)=\phi(v_{11})=1$, $\phi(v_6)=2$ and $\phi(v_{12})=3$, and extend the coloring to $v_4$, $v_5$, $v_9$ and $v_{10}$ by (\refclaim{cl-xext}).
\end{itemize}
\medskip

\noindent{\bf Configuration $\R7.1$.}
If $\phi(v_3)=\phi(v_6)$, then first color the $6$-cycle $v_2v_1v_{10}v_9v_8v_7$ (this is possible, as each of the vertices has at most one colored neighbor), and then color $v_{11}$ and $v_{12}$. Thus, assume that $\phi(v_3)=1$, $\phi(v_6)=2$ and $\phi(v_{12})=3$.  Color the $5$-cycle $v_1v_{11}v_8v_9v_{10}$ (this is possible, as $\phi(x_1)\neq \phi(x_9)$).  Note that in this coloring, $\phi(v_1)\neq 2$ or $\phi(v_8)\neq 1$, as $\phi(v_{11})\neq\phi(v_{12})=3$.  Therefore, the coloring extends to $v_2$ and $v_7$ by (\refclaim{cl-xext}).
\medskip

\noindent{\bf Configuration $\R7.2$.}
The reduction ensures that $\phi(x_1)\neq \phi(x_3)$, say $\phi(x_1)=1$ and $\phi(x_3)=2$.  Also, by symmetry, we may assume that $c=\phi(x_2)\neq 1$. Suppose first that $\phi(v_8)\neq 1$.  Then try coloring $v_{11}$ and $v_3$ by $1$ and $v_1$ by $c$. By (\refclaim{cl-xext}), this coloring extends unless $\phi(v_9)=1$ and $\phi(v_5)=c$.  If $\phi(v_6)\neq 2$, then set the color of $v_3$ to $3$, instead, and observe that the coloring extends.  Otherwise, $\phi(v_6)=2$, and set $\phi(v_{12})=\phi(v_2)=1$, $\phi(v_3)=3$, and color $v_{11}$ and $v_1$.  The coloring extends to
$v_{10}$ and $v_4$ by (\refclaim{cl-xext}).

Therefore, we may assume that $\phi(v_8)=1$.  Suppose that $\phi(v_6)\neq c$.  Then try coloring $v_1$ and $v_{12}$ by $c$, $v_{11}$ and $v_2$ by $5-c$ and $v_3$ by $1$.  By (\refclaim{cl-xext}), this coloring extends to $v_4$ and $v_{10}$ unless $\phi(v_5)=c$ and $\phi(v_9)=1$.  In that case, set $\phi(v_2)=1$, $\phi(v_3)=3$, color $v_{12}$, $v_{11}$ and $v_1$ in this order, and extend the coloring to $v_4$ and $v_{10}$ by (\refclaim{cl-xext}).  Thus, we may assume that $\phi(v_6)=c$.

If $c\neq 2$, then set $\phi(v_3)=c$ and color $v_4$, $v_{10}$, $v_1$, $v_2$, $v_{11}$ and $v_{12}$ in this order; hence, assume that $c=2$.  Consider the coloring that assigns $1$ to $v_2$ and $v_{12}$, $3$ to $v_{11}$ and $v_3$ and $2$ to $v_1$.  If this coloring does not extend to $v_4$ and $v_{10}$, then (\refclaim{cl-xext}) implies that $\phi(v_5)=2$ and $\phi(v_9)=3$.  In that case, set $\phi(v_2)=\phi(v_4)=\phi(v_{12})=1$, $\phi(v_{10})=\phi(v_{11})=2$ and $\phi(v_1)=\phi(v_3)=3$.
\medskip

\noindent{\bf Configuration $\R7.3$.}
The reduction ensures that $\phi(x_1)\neq \phi(x_3)$, say $\phi(x_1)=1$ and $\phi(x_3)=2$.  If $\phi(v_8)\neq 1$ and $\phi(v_6)\neq 2$, then color $v_{11}$ by $1$, $v_{12}$ by $2$ and extend the coloring to the $6$-cycle $v_{10}v_1v_2v_3v_4v_5$.

Assume now that $\phi(v_8)=1$ or $\phi(v_6)=2$.  Suppose first that $\phi(v_6)\neq 2$, and thus $\phi(v_8)=1$. Then try setting the color of $v_1$, $v_5$ and $v_{12}$ to $2$ and coloring $v_{11}$ and $v_{10}$.  If $\phi(x_2)=2$ or $\phi(x_4)=2$ or $\phi(x_2)=\phi(x_4)$, then the coloring extends to $v_2$, $v_3$ and $v_4$, thus assume that $\{\phi(x_2),\phi(x_4)\}=\{1,3\}$. If $\phi(v_9)\neq 2$ or $\phi(v_6)\neq 3$, then set $\phi(v_2)=\phi(v_4)=\phi(v_{11})=2$, $\phi(v_1)=3$, color $v_{12}$ and $v_3$ and extend the coloring to $v_5$ and $v_{10}$ by (\refclaim{cl-xext}).  Otherwise, $\phi(v_9)=2$ and $\phi(v_6)=3$ and we set $\phi(v_5)=1$, $\phi(v_1)=\phi(v_4)=\phi(v_{12})=2$, $\phi(v_{10})=\phi(v_{11})=3$, $\phi(v_2)=\phi(x_4)$ and $\phi(v_3)=\phi(x_2)$.

Therefore, it suffices to consider the case that $\phi(v_6)=2$.  If $\phi(x_4)\neq 2$, then set $\phi(v_4)=2$, color the $5$-cycle $v_1v_2v_3v_{12}v_{11}$, and color $v_{10}$ and $v_5$.  So we have $\phi(x_4)=2$.  Suppose that $\phi(x_2)\neq 2$.  Then set $\phi(v_2)=2$ and $\phi(v_1)=3$. If $\phi(v_8)\neq 2$, then color $v_{11}$ by $2$ and color $v_{10}$, $v_5$, $v_4$, $v_3$ and $v_{12}$ in this order. On the other hand, if $\phi(v_8)=2$, then note that $\phi(v_9)\neq 2$, and set $\phi(v_{10})=2$, $\phi(v_3)=\phi(v_5)=\phi(v_{11})=1$ and $\phi(v_4)=\phi(v_{12})=3$.  Thus, we can assume that $\phi(x_2)=2$.

Try setting $\phi(v_2)=\phi(v_4)=\phi(v_{12})=1$ and $\phi(v_3)=\phi(v_5)=3$. If $\phi(v_9)\neq 1$, then set $\phi(v_{10})=1$ and color $v_{11}$ and $v_1$; thus assume that $\phi(v_9)=1$.  If $\phi(v_8)\neq 2$, then set $\phi(v_{10})=\phi(v_{11})=2$ and $\phi(v_1)=3$.

Finally, consider the case that $\phi(v_9)=1$ and $\phi(v_8)=2$.  Then, we set $\phi(v_3)=\phi(v_5)=\phi(v_{11})=1$, $\phi(v_1)=2$ and $\phi(v_2)=\phi(v_4)=\phi(v_{10})=\phi(v_{12})=3$.
\medskip

\noindent{\bf Configuration $\R7.4$.}
The reduction ensures that $\phi(x_3)\neq \phi(v_6)$, say $\phi(v_6)=1$ and $\phi(x_3)=2$.  Suppose first that $\phi(v_8)\neq \phi(v_{10})$. If $\phi(v_{10})\neq 2$, then let $\phi(v_{12})=2$, $\phi(v_{11})=\phi(v_{10})$ and extend the coloring to the $5$-cycle $v_1v_2v_3v_4v_5$; thus assume that $\phi(v_{10})=2$.  If $\phi(x_2)\neq 2$, then set $\phi(v_2)=2$, $\phi(v_3)=1$, and color $v_4$, $v_5$, $v_1$, $v_{11}$ and $v_{12}$ in this order. If $\phi(x_2)=2$, then set $\phi(v_1)=\phi(v_3)=1$, $\phi(v_2)=3$, and color $v_{11}$, $v_{12}$, $v_4$ and $v_5$, in this order.

Therefore, assume that $\phi(v_8)=\phi(v_{10})=c$.  If $c=2$, then color $v_{12}$ by $2$, extend the coloring to the $5$-cycle $v_1\ldots v_5$, and color $v_{11}$.  If $c=3$, then set $\phi(v_1)=\phi(v_3)=1$, $\phi(v_{11})=2$, $\phi(v_{12})=3$, and color $v_2$, $v_4$ and $v_5$ in this order. Thus, assume that $c=1$.
Try setting $\phi(v_1)=\phi(v_{12})=2$ and $\phi(v_{11})=\phi(v_5)=3$.  If $\phi(x_4)\neq 2$, then set $\phi(v_4)=2$ and color $v_2$ and $v_3$.  If $\phi(x_4)=2$ and $\phi(x_2)\neq 1$, then set $\phi(v_2)=\phi(v_4)=1$ and $\phi(v_3)=3$.

Finally, consider the case that $\phi(x_2)=1$ and $\phi(x_4)=2$.  Then, set $\phi(v_3)=1$, $\phi(v_2)=\phi(v_5)=\phi(v_{11})=2$ and $\phi(v_1)=\phi(v_4)=\phi(v_{12})=3$.
\qed

\section{Discharging}\label{sec-disch}

Let $G$ be a graph in a surface $\Sigma$ with rings $\cal R$.
A face is {\em open $2$-cell} if it is homeomorphic to an open disk.
A face is {\em closed $2$-cell} if it is open $2$-cell and
bounded by a cycle.  A face $f$ is {\em semi-closed $2$-cell} if it is open $2$-cell,
and if a vertex $v$ appears twice in the boundary walk of $f$, then $v$ is the main vertex
of a vertex-like ring $R$ and the edges of $R$ form part of the boundary walk of $f$.
A face $f$ is {\em omnipresent} if it is not open $2$-cell and
each of its boundary walks is a cycle bounding
a closed disk $\Delta\subseteq \widehat{\Sigma}\setminus f$ containing exactly one ring.
We say that $G$ has an {\em internal $2$-cut} if there exist
sets $A,B\subseteq V(G)$ such that $A\cup B=V(G)$, $|A\cap B|=2$,
$A-B\ne\emptyset\ne B-A$, $A$ includes all vertices of $\cal R$,
and no edge of $G$ has one end in $A-B$ and the other in $B-A$.

We wish to consider the following conditions that the triple
 $(G,\Sigma,{\cal R})$ may or may not satisfy:
{
\myitem{(I0)} every internal vertex of $G$ has degree at least three,
\myitem{(I1)} $G$ has no even cycle consisting of internal vertices of
degree three,
\myitem{(I2)} $G$ has no cycle $C$ consisting of internal vertices of degree
three, and two distinct adjacent vertices $u,v\in V(G)-V(C)$ such that both
$u$ and $v$ have a neighbor in $C$,
\myitem{(I3)} every face of $G$ is semi-closed $2$-cell and has length at least $5$,
\myitem{(I4)} if a path of length at most two has both ends in $\cal R$,
then it is a subgraph of $\cal R$,
\myitem{(I5)} no two vertices of degree two in $G$ are adjacent, unless they belong to a vertex-like ring,
\myitem{(I6)} if $\Sigma$ is the sphere and $|{\cal R}|=1$, or if
$G$ has an omnipresent face, then $G$ does not contain an internal $2$-cut,
\myitem{(I7)} the distance between every two distinct members of $\cal R$
is at least four,
\myitem{(I8)} every cycle in $G$ that does not separate the surface has length at least seven,
\myitem{(I9)} if a cycle $C$ of length at most $9$ in $G$ bounds an open disk $\Delta$ in $\Sigma$,
then $\Delta$ is a face, a union of a $5$-face and a $(|C|-5)$-face, or $C$ is a $9$-cycle and $\Delta$ consists of three $5$-faces
intersecting in a vertex of degree three.
}
\bigskip

Let $G$ be a graph in a surface $\Sigma$ with rings $\cal R$ satisfying (I3).
We say that a good configuration $\gamma$ {\em appears} in 
$(G,{\cal R})$ if it faintly appears and the following conditions hold:

{
\myitem{$\bullet$} ${\cal I}_\gamma$ is disjoint from vertex-like rings,
\myitem{$\bullet$} if $\gamma$ is isomorphic to $\R3$, then either ${\cal I}$
contains a vertex of ${\cal R}$ or there exists a vertex $v\in{\cal I}$ such that $v$ and all its neighbors are internal,
\myitem{$\bullet$} if $\gamma$ is isomorphic to $\R4$, then 
the vertex that corresponds to 
$v_2$ is internal and has degree at least $4$, and neither $x_4$ nor $x_5$ belongs to a vertex-like ring,
\myitem{$\bullet$} if $\gamma$ is isomorphic to $\R5$, then
$v_4$ is an internal vertex and
the face whose boundary contains the path corresponding to $v_6v_7v_8$ has 
length at least seven,
\myitem{$\bullet$} if $\gamma$ is isomorphic to $\R6$ or $\R6.1$, then
both vertices of ${\cal A}_\gamma$ are internal, and all neighbors of at least
one of them are internal,
\myitem{$\bullet$} if $\gamma$ is isomorphic to one of 
$\R7,\R7.1,\R7.2,\R7.3,\R7.4$, then
all vertices in ${\cal A}_\gamma\cup{\cal I}_\gamma$ and all their neighbors
are internal, and
\myitem{$\bullet$} if $\gamma$ is isomorphic to $\R7$,
then the vertex corresponding to $x_8$ and all its neighbors are internal.
}

Let $G$ be a graph in a surface $\Sigma$ with rings $\cal R$, 
and let $M$ be a subgraph of $G$ with no isolated vertices.
We define the {\em initial charge} of the triple $(G,\Sigma,{\cal R})$
as follows.
Every face $f$ gets charge $|f|-4$.  A ring vertex of degree two gets
charge $-1$ if it belongs to $M$ and $-1/3$ otherwise, a ring vertex of degree $d\ge 3$ gets charge $d-3$, and all
internal vertices of degree $d$ get charge $d-4$.  
Finally, we increase the charge of each face incident with an edge of $M$
by $5/3$ and each ring vertex of degree two belonging to $M$ by $2/3$.

\begin{lemma}
\mylabel{lem:initcharge}
Let $G$ be a graph in a surface $\Sigma$ with rings $\cal R$,
let $g$ be the Euler genus of $\Sigma$,
let $M$ be a subgraph of $G$ with no isolated vertices, and
let $n_2$ be the number of ring vertices of degree two that do not belong to $M$.
Then the sum of initial charges of all vertices and faces of $G$
is at most 
$4g+4|{\cal R}|+2n_2/3+10|E(M)|/3-8$.
\end{lemma}
\begin{proof}
By Euler's formula, $|E(G)|\le |V(G)|+|F(G)|+|{\cal R}|+g-2$.
Let $n_r$ denote the number of ring vertices.
Note that in the last step of the definition of the initial charge,
we increased the sum of charges by at most $10|E(M)|/3$, since
if $v$ is a ring vertex of degree two belonging to $M$, then
an edge of $M$ incident with $v$ is also incident with only one face of $G$.
The sum of the initial charges of all vertices and faces
is at most
\begin{align*}
&\sum_{v\in V(G)} (\deg(v)-4) + n_r + 2n_2/3 + \sum_{f\in F(G)} (|f|-4)+10|E(M)|/3\\
&=(2|E(G)|-4|V(G)|)+n_r+2n_2/3 + (2|E(G)|-4|F(G)|-n_r)+10|E(M)|/3\\
&=4(|E(G)|-|V(G)|-|F(G)|)+2n_2/3+10|E(M)|/3\\
&\le 4g+4|{\cal R}|+2n_2/3+10|E(M)|/3-8,
\end{align*}
as desired.
\end{proof}

A $5$-face $f$ is {\em $k$-dangerous} if $f$ is not incident with 
an edge of $M$ and
$f$ is incident with exactly $k$ internal vertices of degree three.
Let $f_1=uvawb$ be a $4$-dangerous face, where $w$ is the unique
incident vertex that is not internal of degree three.
Let $f_2$ be the face incident with $uv$ distinct from $f_1$.
We say that $f_2$ is {\em linked to $f_1$ (through the edge $uv$)}.
Let $xy$ be an edge such that $y$ has degree three, and let $g_1$, $g_2$, $g_3$
be the faces incident with $y$ such that $xy$ is incident with $g_1$ and $g_2$.
Then the face $g_3$ is {\em opposite} to $x$.
A $4$-dangerous face $f$ is {\em extremely $4$-dangerous} if it is neither incident with
a vertex of ${\cal R}$ nor opposite to the main vertex of a vertex-like ring.

Let us apply the following {\em primary discharging rules}, 
resulting in the {\em primary charge}:

\begin{description}
\item[Rule~1:] Every face sends $1/3$ to each incident ring vertex 
of degree two and each incident internal vertex of degree three.
\item[Rule~2:] If $uvw$ is a subpath of a ring, then $v$ sends $1/3$ to each 
face incident with $v$ other than the two faces incident with $uv$ and $uw$.
Additionally, if $v$ is the main vertex of a vertex-like ring, then $v$ sends
$1/3$ to each opposite face and receives $2/3$ from the face incident with the ring.
\item[Rule~3:] Let $f$ be a face linked to an extremely $4$-dangerous face $f'$ through an edge $uv$.
If $f$ has length at least $6$, or $f$ is incident with an edge of $M$, then $f$ sends $1/3$ to $f'$ across the edge $uv$.
\item[Rule~4:] Let $v_1v_2v_3v_4$ be a subwalk of the boundary walk of a
face $f'$ of length at least seven,
such that $f'$ is linked to extremely $4$-dangerous faces through both 
$v_1v_2$ and $v_3v_4$.  
Let $f$ be the other face incident with the edge $v_2v_3$.
If $f$ has length at least six, 
then $f$ sends $1/9$ to $f'$ across the edge $v_2v_3$.
\end{description}

\begin{lemma}
\mylabel{lem:primaryvertex}
Let $G$ be a graph in a surface $\Sigma$ with rings $\cal R$
satisfying {\rm (I0)} and {\rm (I3)}
and let $M$ be a subgraph of $G$.
Then the primary charge of each vertex is non-negative,
and the primary charge of a ring vertex of degree $d\ge 4$ is at least $(d-2)/3$.
Moreover, the primary charge of an internal vertex of degree $d\ge4$
is exactly $d-4$.
\end{lemma}
\begin{proof}
By Rule~1, the internal vertices of degree three have primary charge $0$.
The charge of internal vertices of degree
$d\ge 4$ is unchanged, i.e., $d-4\ge 0$.
Consider now a ring vertex $v$ of degree $d$.
If $d=2$, then the initial charge of $v$ is $-1/3$ and $v$ receives $1/3$ by Rule~1.  
Observe that $v$ sends nothing by Rule~2, thus the primary charge of $v$ is $0$.
If $d\ge 3$, then $v$ sends charge by Rule~2 to $d-3$ incident faces.
Furthermore, if $v$ is the main vertex of a vertex-like ring, then $v$ sends $1/3$
to at most $d-2$ opposite faces and receives $2/3$ from the face incident with the ring.
Hence, the primary charge of $v$ is at least $d-3-\max((d-3)/3,(2d-7)/3)$,
which is non-negative, and at least $(d-2)/3$ for $d\ge 4$ as desired.
\end{proof}

Let us now estimate the primary charge of faces.
A subgraph $M\subseteq G$ \emph{captures $(\le\!4)$-cycles} if $M$ contains all
cycles of $G$ of length at most $4$ and furthermore, $M$ is either null or has minimum
degree at least two.

\begin{lemma}
\mylabel{lemma-primary}
Let $G$ be a graph in a surface $\Sigma$ with rings $\cal R$ 
satisfying {\rm (I0)}, {\rm(I1)}, {\rm(I3)}, {\rm(I4)}, {\rm(I5)} and {\rm(I7)},
let $M$ be a subgraph of $G$ that captures $(\le\!4)$-cycles
and assume that if a configuration isomorphic to one of $\R1$, $\R2$, \ldots, $\R5$
appears in $G$, then it touches $M$.
If $f$ is a face of $G$, then the primary charge of $f$ is non-negative.  
Furthermore, if the primary charge of $f$ is zero, then $f$ has length
exactly five, it is not incident with an edge of $M$, and
{
\myitem{(a)} $f$ is $3$-dangerous, or
\myitem{(b)} $f$ is incident with a ring vertex, or
\myitem{(c)} $f$ is $4$-dangerous and a face of length at least $6$ is linked to $f$, or
\myitem{(d)} $f$ is $4$-dangerous, the face $h$ linked to $f$ has length 
five and $h$ is incident with an edge of $M$, or
\myitem{(e)} $f$ is $4$-dangerous and is opposite to the main vertex of a vertex-like ring.
}

\noindent
Otherwise, the primary charge of $f$ is least $2/9$, and if $|f|\ge 8$,
then the primary charge of $f$ is at least $5|f|/9-4$.  Also,
if $f$ is a $6$-face incident with a ring vertex of degree two, then $f$ has primary charge at least $2/3$.
\end{lemma}
\begin{proof}
Suppose first that $f$ has length exactly five.  The face $f$ may send charge by Rules~1 and 3.
Let us consider the case that $f$ is incident with an edge of $M$.
If $f$ sends charge across an edge $uv$ by Rule~3 to a face $f'$, then both $u$ and $v$ have degree three
and no edge of $f'$ belongs to $M$.  Since $M$ has minimum degree at least two, it follows that no edge
incident with $u$ or $v$ belongs to $M$; hence $f$ sends charge by Rule~3 to at most two faces.
The primary charge of $f$ is at least $1+5/3-5/3-2/3=1/3>2/9$.

Therefore, we may assume that $f$ is not incident with any edge of $M$, and in particular, $f$ does not share
an edge with any cycle of length at most $4$.  Also, $f$ sends charge only by Rule~1.
Let us distinguish several cases according to the number of internal vertices 
of degree three incident with $f$.

\medskip

\noindent $\bullet$
{\em All vertices incident with $f$ are internal and have degree three.}
Then $f$ and its incident vertices form a configuration isomorphic to $\R1$
that appears in $G$, which is a contradiction.

\medskip

\noindent  $\bullet$
{\em The face $f$ is incident with exactly four 
internal vertices of degree three.} 
Let $f=v_1v_2v_3v_4v_5$ and suppose that all these vertices except for 
$v_2$ are internal and have degree three.  If $v_2$ is not internal,
then $v_2$ has degree at least four, since
$v_1$ and $v_3$ are internal vertices.  The charge of $f$ after applying Rule~1 is $-1/3$.

The face $f$ is incident with no edge of $M$, hence $f$ is $4$-dangerous.  
If $v_2$ belongs to a ring, then $f$ receives $1/3$ by Rule~2, 
making its charge zero, and hence $f$ satisfies (b).
Thus we may assume that $v_2$ is internal and of degree at least $4$.
Similarly, if $f$ is opposite to the main vertex of a vertex-like ring, then $f$ receives $1/3$ by Rule~2
and $f$ satisfies (e), hence it suffices to consider the case that $f$ is extremely $4$-dangerous.

If the face $h$ with that $f$ shares the edge $v_4v_5$ has length five,
then the faces $f$ and $h$ form an imprint of $\R4$ ($v_2$ is distinct from
the vertices incident with $h$, since $f$ does not share an edge with a cycle of
length at most $4$), and a configuration isomorphic to $\R4$ appears in $G$.
By hypothesis the face $h$ is incident with an edge of $M$.

We conclude that $h$ either has length at least $6$ or is incident with an edge of $M$.
In both cases, $h$ sends $1/3$ to $f$ by Rule~3.  Thus the primary charge
of $f$ is zero, and $f$ satisfies (c) or (d).

\medskip

\noindent  $\bullet$
{\em The face $f$ is incident with exactly three internal vertices 
of degree three.} 
In this case $f$ sends $1/3$ to each of the three incident internal
vertices of degree three by Rule~1, making its charge zero.
(The face $f$ is not incident with a ring vertex of degree 
two, since both neighbors of such a vertex belong to $\cal R$).
Since $f$ does not share an edge with $M$, $f$ is $3$-dangerous and satisfies (a).

\medskip

\noindent  $\bullet$
 {\em The face $f$ is incident with exactly two internal vertices 
of degree three.} 
Then $f$ sends $1/3$ to each of them, and at most $1/3$ to a ring vertex 
of degree two by Rule~1,
making its charge non-negative.  
Furthermore, if the charge is zero, then $f$ satisfies (b); otherwise
the charge is at least $1/3$, as desired.

\medskip

\noindent  $\bullet$
 {\em The face $f$ is incident with at most one internal vertex of degree 
three.} 
Then $f$ sends at most $2/3$ by Rule~1 and (I5), 
and its primary charge is at least $1/3$, as desired.

\bigskip

Thus we have proved the lemma when $f$ has length five.
Let us now consider the case that $f$ has length six, and let $f=v_1v_2v_3v_4v_5v_6$.
By (I1) not all vertices incident with $f$ are internal and of degree three.
Thus $f$ sends at most $5/3$ by Rule~1 and at most $4/3$ by Rules~3 and 4;
furthermore, if $f$ sends $2/3$ by Rule~2 (i.e., a vertex-like ring forms part of the boundary of $f$),
then $f$ sends at most $4/3$ by Rule~1 and at most $1/3$ by Rules~3 and 4.
If $f$ is incident with an edge of $M$, then its primary charge is at least
$2+5/3-5/3-4/3=2/3$, as desired, and so we may assume that $f$ is incident
with no edge of $M$. Since $M$ captures $(\le\!4)$-cycles, it follows that no edge of $f$
is incident with a vertex-like ring.

If, say, $v_1$ is the main vertex of a vertex-like ring, then {\rm (I7)} implies that all other
vertices incident with $f$ are internal.  Also, observe that for each of the edges $v_1v_2$, $v_1v_6$, $v_2v_3$ and $v_5v_6$,
either not both ends of the edge are internal vertices of degree three, or the edge is not incident with an extremely $4$-dangerous face;
hence, $f$ sends at most $2/3$ by Rule~3 and nothing by Rule~4.
Furthermore, $f$ receives $1/3$ from $v_1$ by Rule~2, and thus the primary charge of $f$ is at least $2-5/3-2/3+1/3=0$.
If $f$ sends less than $5/3$ by Rule~1 or less than $2/3$ by Rule~3, then the primary charge is at least $1/3$, as desired.
Otherwise, $f$ forms an appearance of $\gamma=\R3$, with ${\cal I}_\gamma=\{v_2,v_4,v_6\}$,
contradicting the hypothesis of the lemma.  Therefore, no vertex incident with $f$ is the main vertex of a vertex-like ring.

Suppose that $f$ sends charge across $v_2v_3$ by Rule~3 or 4.  It follows that $v_2$ and $v_3$ are internal and of degree three.
Let $x_2$ be the neighbor of $v_2$ other than $v_1$ and $v_3$, and let $x_3$ be defined analogously.
Then both $x_2$ and $x_3$ are internal vertices of degree three.  
If $v_1$ and $v_5$ both belong to $\cal R$,
then by {\rm (I4)} $v_6$ is a vertex of degree two, and by {\rm(I4)} and
{\rm (I5)} $v_4$ is an internal vertex, implying that $\gamma=\R3$ appears in $G$ (with ${\cal I}_\gamma=\{v_2,v_4,v_6\}$).
This contradicts the hypothesis; hence, assume that at least one of $v_1$ and $v_5$ is internal, and symmetrically, at least
one of $v_4$ and $v_6$ is internal.
If both $v_1$ and $v_5$ are internal, then then $\gamma=\R3$ appears in $G$ with ${\cal I}_\gamma=\{v_2,v_4,v_6\}$.
And if exactly one of $v_1$ and $v_5$ belongs to $\cal R$, then $\gamma=\R3$ appears in $G$
with ${\cal I}_\gamma=\{v_1,v_3,v_5\}$.  This is a contradiction, showing that $f$ does not send charge across $v_2v_3$ by Rule~3 or 4.

By symmetry, $f$ does not send charge using Rules~3 or 4 at all, and thus its primary charge
is at least $2-5/3=1/3$.  Furthermore, if some vertex incident with $f$, say $v_2$, has degree two and belongs to a ring $R$,
then by (I5), $v_1$ and $v_3$ belong to $R$ and have degree at least three, and thus $f$ sends at most $4/3$ by Rule~1,
and the primary charge of $f$ is at least $2/3$.  This completes the case $|f|=6$.

\bigskip

Finally, we consider the case that $|f|\ge 7$.
Let us estimate the amount of charge sent from $f$ and received by 
$f$ using Rules~3 and 4.  If $v_1v_2v_3v_4$
is a subwalk of the boundary walk of $f$ and $f$ sends $1/3$ across 
$v_2v_3$ by Rule~3, then assign $1/9$ of this charge
to each of $v_1v_2$, $v_2v_3$ and $v_3v_4$.  If $f$ sends $1/9$ across $v_2v_3$ by Rule~4, then add $1/9$ to the charge
assigned to $v_2v_3$; if $f$ receives $1/9$ across $v_2v_3$, then remove $1/9$ from the charge assigned to $v_2v_3$.
We claim that each edge has at most $1/9$ assigned to it, and hence that $f$ loses at most $|f|/9$ by Rules~3 and 4.

Suppose for a contradiction that more than $1/9$ is assigned to the edge $v_2v_3$.  By symmetry, we
can assume that $f$ sends charge by Rule~3 to the face $f_{12}$ across $v_1v_2$.
Let $f_{23}\neq f$ be the face incident with the edge $v_2v_3$.  
If $f$ sends charge across $v_2v_3$ by Rule~3, then
the faces $f_{12}$ and $f_{23}$ form an appearance of 
a configuration isomorphic to $\R5$.
It follows that $f_{12}$ or $f_{23}$ is incident with an edge of $M$.  
This is a contradiction, because Rule~3 sends charge to $4$-dangerous faces
only.  Furthermore, $f$ does not send charge across $v_2v_3$ by Rule~4, because
$f$ is linked to $f_{12}$ through $v_1v_2$.

Since more than $1/9$ is assigned to $v_2v_3$, it follows that $f$ sends charge across
$v_3v_4$ by Rule~3 and does not receive charge by Rule~4 across $v_2v_3$.  Therefore,
$f_{23}$ has length five and $f_{12}$ and $f_{23}$ form an appearance of a configuration
isomorphic to $\R5$ as before.  Since $f_{12}$ is $4$-dangerous, some edge of $M$ is incident with $f_{23}$
but not with $f_{12}$.
Since all neighbors of $v_2$ and $v_3$ have degree three and $M$ has
minimum degree at least two, it follows that some edge of $M$ is incident
with the face $f_{34}\neq f$ that is incident with $v_3v_4$.  
This is a contradiction, because $f$ sends charge to $f_{34}$ by Rule~3.

We can now bound the primary charge of $f$.  If $f$ has length at least eight,
then $f$ sends at most $|f|/3$ by Rule~1 and at most $|f|/9$ by Rules~3 and 4 (and any charge sent by Rule~2
is dominated by the charge received due to sharing an edge with $M$); thus its primary charge
is at least $|f|-4-|f|/3-|f|/9=5|f|/9-4>2/9$, as desired.

Finally, assume that $f$ has length exactly seven.
If $f$ is incident with an edge of $M$,
then $f$ sends at most $7/3$ by Rule~1, 
making the primary charge of $f$ at least $3+5/3-7/3-7/9=14/9$.  
If $f$ is incident with no edge of $M$, then
$f$ and its incident vertices do not form an
appearance of a configuration isomorphic to $\R2$, and that
in turn implies
that $f$ is incident with no more than six internal vertices of degree three.
Thus $f$ sends at most $2$ by Rule~1, and hence
the primary charge of $f$ is at least $3-2-7/9=2/9$, as desired.
\end{proof}

\begin{figure}
\begin{center}
\epsfbox{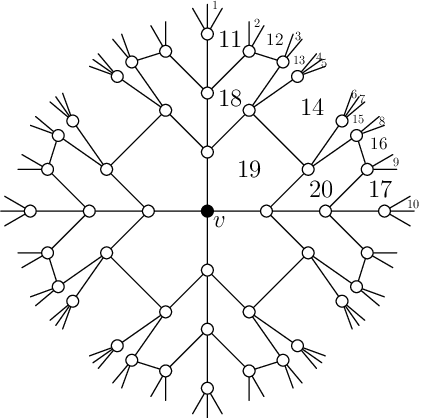}
\end{center}
\caption{The maximum number of $3$-reachable faces (case $d=4$).}
\label{fig-20d}
\end{figure}

We now modify the primary charges using three additional rules
into what we will call ``final charges".
A vertex is {\em safe} if its degree is at least five, 
or if it belongs to $\cal R$, or if it is incident with a face with 
strictly positive primary charge.
A face $f$ is {\em $k$-reachable} from a vertex $v$ if there exists a path $P$
of length at most $k$ ($P$ may have length zero), joining $v$ to a vertex incident with $f$, such that no vertex of
$P\backslash v$ is safe.  In particular, every vertex of $P\backslash v$ is internal and has degree at most four,
and all faces incident with them have length $5$, which implies that the number of faces that are
$3$-reachable from a vertex of degree $d$ is bounded by $20d$ (see Figure~\ref{fig-20d} demonstrating the worst case).
Furthermore, if $v$ is a ring vertex or an internal vertex incident to
a face $f$ with strictly positive primary charge, then two of the neighbors of $v$ are safe, and we conclude that
at most $20(d-3)+26$ faces distinct from $f$ are $3$-reachable from~$v$.

Let $\epsilon>0$ be a real number, to be specified later.
Starting from the primary charges we now apply the following three
rules, resulting in the {\em final charge}:

\begin{description}
\item[Rule~5:] 
The charge of each ring vertex of degree three is increased by $26\epsilon$.
\item[Rule~6:] Each face of strictly positive primary charge sends
$46\epsilon$ units of charge to each incident vertex.
\item[Rule~7:] If $v$ is a safe vertex of degree at least three, then
$v$ sends a charge of $\epsilon$ to each face of zero primary charge that is $3$-reachable from $v$.
\end{description}

\begin{lemma}
\mylabel{lem:fincharge}
Let $G$ be a graph in a surface $\Sigma$ with rings $\cal R$,
let $g$ be the Euler genus of $\Sigma$,
let $M$ be a subgraph of $G$ that captures $(\le\!4)$-cycles,
let $n_2$ be the number of ring vertices of degree two not belonging to $M$,
let $n_3$ be the number of ring vertices of degree three,
let $\epsilon>0$, and let $M$ be a subgraph of $G$.
Then the sum of final charges of all vertices and faces of $G$
is at most
$4g+4|{\cal R}|+26\epsilon n_3+2n_2/3+10|E(M)|/3-8$.
\end{lemma}

\proof
This follows directly from Lemma~\ref{lem:initcharge} and the description of
the discharging rules.~\qed
\bigskip

\begin{lemma}
\mylabel{lem:finalvertex}
Let $G$, $\Sigma$, ${\cal R}$, and $M$ be as in Lemma~{\rm\ref{lemma-primary}},
and let $\epsilon\le 1/180$.
Then the final charge of every vertex is non-negative and
the final charge of every ring vertex of degree $d\ge 4$
is at least $(1/3-20\epsilon)(d-2)-26\epsilon$.
\end{lemma}

\proof
Let $v$ be a vertex of $G$ of degree $d$.
Lemma~\ref{lem:primaryvertex} tells us that the primary charge of 
$v$ is non-negative.
If $v$ is safe, then it sends at most $20\epsilon d$ units of charge
by Rule~7; otherwise it sends nothing using Rules~5--7.
Assume first that $v$ is an internal vertex.
If $d\ge5$, then the primary charge of $v$ is $d-4$, and its final
charge is at least $d-4-20\epsilon d$, which is non-negative by
the choice of $\epsilon$.
If $d\le4$ and $v$ is not incident with a face of positive primary
charge, then its final charge is the same as its primary charge,
and so the conclusion follows from Lemma~\ref{lem:primaryvertex}.
If $d\le4$ and $v$ is incident with a face of positive primary charge,
then it receives at least $46\epsilon$ units of charge using Rule~6
and sends at most $46\epsilon$ units using Rule~7.
Thus $v$ has non-negative final charge.

Let us now assume that $v$ is a ring vertex.
If $d=2$, then $v$ sends no charge by Rules~5--7 and its final charge is zero.
If $d=3$, then $v$ receives $26\epsilon$ units using Rule~5,
and sends at most $26\epsilon$ units using Rule~7.
Finally, if $d\ge4$, then $v$ has primary charge at least $(d-2)/3$
by Lemma~\ref{lem:primaryvertex}, and it sends at most $20(d-3)\epsilon+26\epsilon$
units of charge, and hence its final charge is at least
$(1/3-20\epsilon)(d-3)+1/3-26\epsilon$, which is non-negative by the choice of $\epsilon$.~\qed
\bigskip

\begin{lemma}
\mylabel{lem:finalbigface}
Let $G$, $\Sigma$, ${\cal R}$, and $M$ be as in Lemma~{\rm\ref{lemma-primary}},
and let $\epsilon>0$ be arbitrary.
Then the final charge of every face of length six or seven is at least $2/9-322\epsilon$,
and the final charge of every face of length $l\ge8$ is
at least $(5/9-46\epsilon)l-4$.
\end{lemma}

\proof
Lemma~\ref{lemma-primary} gives a lower bound on the primary charge of a face 
$f$, and $f$ sends at most $46\epsilon|f|$ units of charge 
using Rule~6.~\qed
\bigskip

\begin{lemma}
\mylabel{lem:safereach}
Let $G$, $\Sigma$, ${\cal R}$, and $M$ be as in Lemma~{\rm\ref{lemma-primary}}, satisfying additionally {\rm (I8)},
and assume that if a configuration isomorphic to one of 
$\R1$, $\R2, \ldots, \R6$ or $\R7$
appears in $G$, then it touches $M$.
Then every face of zero primary charge is $3$-reachable from some safe vertex.
\end{lemma}

\proof 
Let $f$ be a face of zero primary charge.
Lemma~\ref{lemma-primary} implies that $f$ is a $5$-face, and unless $f$ is $1$-reachable from a safe vertex,
we have that $f$ is $3$-dangerous and all vertices incident with $f$ are internal and have degree at most four.
Let $f=w_1w_2w_3w_4w_5$, and suppose first that $w_1$ and $w_5$
have degree four.  In this case, we prove the following stronger claim:
both $w_1$ and $w_5$ are at distance at most two from a safe vertex.

Let $f'$ be the other face incident with the edge $w_1w_5$.
To prove the claim we may assume that no vertex incident with $f$ or $f'$
is safe, for otherwise the claim holds.
Then $f'$ has primary charge zero, because no vertex incident with $f$ is safe.
Since $w_1$ and $w_5$ have degree at least four,
Lemma~\ref{lemma-primary} implies that $f'$ is $3$-dangerous.
Since $f$ and $f'$ have zero primary charge,
they do not share an edge with $M$, and in particular, they do not share an edge
with any cycle of length at most four.
We deduce that the faces $f$ and $f'$ and their incident vertices form a
faint appearance of a configuration isomorphic to $\R6$.
Since $f$ and $f'$ are incident with no edge of $M$, this is not
an appearance; hence either $w_1$ or $w_5$ has a neighbor in $\cal R$, or
the distance from both $w_1$ and $w_5$ to a vertex of $\cal R$
is most two.  In both cases, $w_1$ and $w_5$ are at distance at most two from a safe vertex, as desired.
This concludes the case when $w_1$ and $w_5$ have degree four.

We may therefore assume that $w_1$ and $w_3$ have degree four.
Let $f_1$, $f_2$,  $f_3$, $f_4$ and $f_5$ be the other faces incident with the
edges $w_1w_2$, $w_2w_3$, $w_3w_4$, $w_4w_5$ and $w_5w_1$, respectively.
Similarly as before we may assume that 
$f_1$, $f_2$, $f_3$, $f_4$ and $f_5$ are all $3$-dangerous $5$-faces and
vertices incident with them have degree at most four, for otherwise 
$f$ is $3$-reachable from a safe vertex.
If any of those faces contained two consecutive vertices $x$ and $y$ of degree four, 
then by the previous paragraph, both $x$ and $y$ would be at distance at most
two from a safe vertex, and hence $f$ would be $3$-reachable from such a safe vertex.
We may therefore assume that this is not the case.  Since no cycle of length at most $4$
shares an edge with $f$ or $f_i$ for $1\le i\le 5$, we deduce that the
faces $f,f_1,f_2,f_3,f_4,f_5$ and their 
incident vertices and edges form a faint
appearance of a configuration $\gamma$ isomorphic to $\R7$, unless
$f_3$ and $f_5$ are incident with a common vertex, i.e., unless
$v_4$ is identified with $v_9$, or $v_5$ is identified with $v_{10}$
in the depiction of $\R7$ in Figure~\ref{fig-redu}.  Suppose that say
$v_4=v_9$.  Since this vertex has degree three, we conclude that
$\{v_3,v_5\}\cap \{v_8,v_{10}\}\neq \emptyset$.  As $f$ does not share an edge
with $M$, we have $v_3\neq v_8$, $v_3\neq v_{10}$ and $v_5\neq v_8$.  However, if $v_5=v_{10}$,
then the cycle $v_5v_6v_{12}v_{11}v_1$ 
does not separate the surface, contrary to {\rm (I8)}.

It follows that $\R7$ faintly appears, but does not appear, in $G$.
Thus, using the labeling of the vertices as in Figure~\ref{fig-redu},
one of $x_1$, $x_3$, $x_6$, $x_7$, $x_8$ or one of their neighbors
belongs to $\cal R$.  Therefore, $f$ is $3$-reachable from a safe
vertex, as desired.~\qed
\bigskip

\begin{lemma}
\mylabel{lem:final5face}
Let $G$, $\Sigma$, ${\cal R}$, and $M$ be as in Lemma~{\rm\ref{lem:safereach}},
let $\epsilon\le2/2079$,
and assume that if a good configuration appears in $G$,
then it touches $M$.
Then the final charge of every face of length five is at least 
$\epsilon$.
\end{lemma}

\noindent
Let us remark that $2079=9(5\cdot46+1)$.

\proof
Let $f$ be a face of length five.
If $f$ has positive primary charge, then by Lemma~\ref{lemma-primary}
it has primary charge at least $2/9$. 
It sends $46\epsilon$ units of charge to each incident vertex by Rule~6,
and hence $f$ has final charge at least $2/9-5\cdot 46\epsilon\ge \epsilon$.

We may therefore assume that $f$ has primary charge zero.
By Lemma~\ref{lem:safereach}, $f$ is $3$-reachable from some safe vertex,
and hence has final charge at least $\epsilon$ because of Rule~7, 
as desired.~\qed
\bigskip

Let $s:\{5,6,\ldots\}\to{\mathbb R}$ be a function (that we specify later) satisfying
{
\myitem{(S1)} $s(5)= 2\epsilon$,
\myitem{(S2)} $s(7)\le 4/9-644\epsilon$, and
\myitem{(S3)} $s(l)\le (10/9-92\epsilon)l-8$ for every integer $l\ge8$.
}

\noindent Suppose that we are given such a function and a graph $G$ in $\Sigma$ with rings $\cal R$.
For a face $f$ of $G$, we define $w(f)=s(|f|)$ if $f$ is open $2$-cell,
and $w(f)=|f|$ otherwise.  We define $w(G,{\cal R})$ as the sum of $w(f)$ over all faces
$f$ of $G$.

\begin{lemma}
\mylabel{lem:fincharges}
Let $G$ be a graph in a surface $\Sigma$ with rings $\cal R$
satisfying {\rm (I0)--(I8)},
let $M$ be a subgraph of $G$ that captures $(\le\!4)$-cycles
and assume that if a configuration isomorphic to one of $\R1,\R2,\ldots,\R7$
appears in $G$, then it touches $M$.
Let $\epsilon$ be a real number satisfying $0<\epsilon<2/2079$,
and let $s:\{5,6,\ldots\}\to{\mathbb R}$ be a function satisfying
{\rm (S1)--(S3)}.
Then the final charge of every vertex is non-negative, 
and the final charge of every face $f$ is at least $s(|f|)/2$.
\end{lemma}

\proof
The assertions follow from Lemmas~\ref{lem:finalvertex},
\ref{lem:finalbigface} and~\ref{lem:final5face} using conditions
(S1)--(S3).~\qed

\begin{lemma}
\mylabel{lem:noconfigweight}
Let $G$, $\Sigma$, ${\cal R}$, $M$, $\epsilon$, and $s$ be as in Lemma~{\rm\ref{lem:fincharges}},
and let $g$, $n_2$, and $n_3$ be as in Lemma~\ref{lem:fincharge}.
Then $w(G,{\cal R})\le 8g+8|{\cal R}|+52\epsilon n_3+4n_2/3+20|E(M)|/3-16$.
\end{lemma}

\proof
By Lemma~\ref{lem:fincharges} the quantity $w(G,{\cal R})$ is at most
twice the sum of the final charges of all vertices and faces of $G$,
and hence the lemma follows from Lemma~\ref{lem:fincharge}.~\qed
\bigskip

We need the following refinement of the previous lemma.
\begin{lemma}
\mylabel{lem:newconfigweight}
Let $G$, $\Sigma$, ${\cal R}$, $M$, $\epsilon$, and $s$ be as in 
Lemma~{\rm\ref{lem:fincharges}}, and let $g$, $n_2$, and $n_3$ be as in Lemma~\ref{lem:fincharge}.
Then $w(G,{\cal R})\le 8g+8|{\cal R}|+52\epsilon n_3+4n_2/3+20|E(M)|/3-8b/9-16$,
where $b$ is the number of $6$-faces of $G$ incident with a ring vertex of degree
two, plus the number of ring vertices of degree at least four.
\end{lemma}

\proof
This follows similarly as Lemma~\ref{lem:noconfigweight}, since according to Lemma~\ref{lemma-primary},
each $6$-face incident with a ring vertex of degree two has charge by at least $4/9$ higher than the bound
used to derive Lemma~\ref{lem:noconfigweight}, and since the final charge of a ring vertex of degree at least four
is at least $2/3-46\epsilon > 4/9$.~\qed

\section{Reductions}\label{sec-reductions}

In this section, we argue that subject to a few assumptions, reducing
a good configuration does not create cycles of length at most four.

Let $G$ be a graph in a surface $\Sigma$ with rings $\cal R$,
and let $P$ be a path of length at most four with ends
$u,v\in V(\cal R)$ and otherwise disjoint from $\cal R$.
We say that $P$ is {\em allowable} if
\begin{itemize}\addtolength{\itemsep}{-7pt}
\item $u,v$ belong to the same ring of $\cal R$, say $R$,
\item $P$ has length at least three,
\item there exists a subpath $Q$ of $R$ with ends $u,v$
such that $P\cup Q$ is a cycle of length at most eight
that bounds an open disk $\Delta\subset \Sigma$,
\item if $P$ has length three, then $P\cup Q$ has length five
and $\Delta$ is a face of $G$, and
\item if $P$ has length four, then $\Delta$ includes at
most one edge of $G$, and if it includes such an edge $e$, then
$Q$ has length four and $e$ joins the middle vertex of $P$ to the middle vertex of $Q$.
\end{itemize}

\noindent
We say that $G$ is {\em well-behaved} if every path $P$ of length 
at least one and at most four with ends 
$u,v\in V(\cal R)$ and otherwise disjoint from
$\cal R$ is allowable.

We say that a configuration $\gamma$ {\em strongly appears in $G$} if it both appears and weakly appears in $G$
and
\begin{itemize}\addtolength{\itemsep}{-7pt}
\item if $u,v\in{\cal A}_\gamma$ are distinct, 
then at least one of $u,v$ is internal,
\item if $u,v\in{\cal I}_\gamma$ are distinct,
$u\in V({\cal R})$, and $w\in V({\cal R})$ is a neighbor of $v$, then
$u$ and $w$ are adjacent and $uw, wv\in E(G_{\gamma})$, and
\item if $\gamma$ is isomorphic to $\R7$, then 
the vertices corresponding to $v_2$ and $z$ are distinct,
non-adjacent and have no common neighbor distinct from $v_1$, $v_3$, $x_6$ and $x_7$.
\end{itemize}

\begin{lemma}
\mylabel{lem:strongappear}
Let $G$ be a graph in a surface $\Sigma$ with rings $\cal R$ satisfying {\rm (I0)}, {\rm (I2)} and {\rm (I8)},
and assume that $G$ is well-behaved.
If a configuration isomorphic to one of $\R1$, $\R2$, \ldots, $\R7$ appears in $G$
and no cycle in $G$ of length four or less touches it,
then either a good configuration strongly appears in $G$,
or $\Sigma$ is a disk, ${\cal R}=\{R\}$, $R$ has
length $2s$ for some $s\in \{5,7\}$, 
$V(G)=V(R)\cup V(C)$ for a cycle $C$ of length $s$, 
and each vertex of $C$ is internal
of degree three and has one neighbor in $R$.
\end{lemma}

\proof
Let $\gamma$ be a good configuration appearing in $G$,
such that no cycle in $G$ of length four or less touches $\gamma$.
If possible, we choose $\gamma$ so that it is equal to one of $\R1$, $\R6.1$, $\R7.1$, $\R7.2$, $\R7.3$ or $\R7.4$.
We claim that, possibly after relabeling the vertices of $G_\gamma$,
$\gamma$ strongly appears in $G$.
To prove that we first notice that the first condition of weak
appearance holds by hypothesis and (I8)---if $x_3=x_7$, then $x_3v_3v_{12}v_6v_7$ is a $5$-cycle
separating $x_1$ from $x_6$.  The third condition is implied by appearance.
The second condition of weak appearance follows from
our choice of $\gamma$ and the fact that no cycle of length
at most four touches $\gamma$.  For example, if $\gamma$ is $\R7$, then $v_2$ and $v_7$ are not adjacent,
because $\R7.1$ does not appear in $G$ by the choice of $\gamma$.
Additionally, when $\gamma$ is $\R5$, we use (I2) to show that $v_1$ is not adjacent to $v_5$.

It remains to prove that $\gamma$ satisfies the conditions
of strong appearance.  Let us discuss the configurations separately.
If $\gamma$ is $\R1$ or $\R2$, it suffices to show that we can choose the labels
of the vertices of $\gamma$ so that $x_1$ is internal.  If that is not possible, then
each vertex of $\gamma$ is adjacent to a vertex belonging to $\cal R$.
Since $G$ is well-behaved it follows that there exists a ring $R\in{\cal R}$
that satisfies the conclusion of the lemma for $s=5$ if $\gamma$ is $\R1$
and for $s=7$ when $\gamma$ is $\R2$.

If $\gamma$ is $\R3$, we only need to prove the second condition of strong appearance.
Suppose that say $v_3\in V({\cal R})$ and $v_5$ has a neighbor $x_5$
in ${\cal R}$ other than $v_4$.  Since $G$ is well-behaved, $v_4$ is an internal vertex
and $v_3v_4v_5x_5$ together with a path in ${\cal R}$ bound a $5$-face, implying that $v_4$
has degree two.  This contradicts {\rm (I0)}.  

If $\gamma$ is $\R4$, then note that the path $x_1v_1v_2v_3x_3$ 
is not allowable, since by the definition of appearance,
$v_2$ has degree at least four.   
Therefore, at least one of $x_1$ and $x_3$ is internal, and $\gamma$ strongly appears.

If $\gamma$ is $\R5$, we need to prove the first and the second condition of strong appearance.
For the first one, observe that the path $v_2v_1v_8x_8$ is not allowable,
since $v_1$ has degree at least three.  For the second condition,
since $\gamma$ appears in $G$, we have that $v_4$ is internal; thus
it suffices to consider the case that $x_6$ and a neighbor $x_4$ of $v_4$ belongs to $\cal R$.
Since $v_3v_4v_5v_6v_7$ is not an appearance of $\R1$ in $G$, $v_4$ has degree at least four, and
thus the paths $v_2v_3v_4x_4$ and $x_4v_4v_5v_6x_6$ cannot both be allowable.  It follows that
$v_2$ is internal, and similarly all neighbors of $v_2$ are internal.  However, then we can
relabel the vertices of $\gamma$, switching $v_2$ with $v_4$, $v_6$ with $v_8$, etc., and obtain
a strong appearance of $\R5$ in $G$.

For the configurations $\R6$, \ldots, $\R7.4$, the first two conditions follow from the definition
of appearance.  Therefore, suppose that $\gamma$ is $\R7$ and let us now consider the last condition
in the definition of strong appearance.  Again, we we use symmetry: if the condition does not hold for $\gamma$ we
swap $v_1$ and $v_3$, $v_6$ and $v_8$, and so on.
The vertex $v_2$ cannot be equal to or adjacent to both $z$ and $z_1$,
since $v_2\neq x_7$ (otherwise, $\R7.1$ would appear in $G$), $x_7$ has degree at least three and no cycle of length at most four touches $\gamma$.  Unless the condition holds,
we can assume that $z_1\neq v_2$, $z_1$ is not adjacent to $v_2$ and that $z_1$ and $v_2$ have a common neighbor $x_2$ distinct
from $v_1$, $v_3$, $x_7$ and $x_8$.  Since no cycle of length at most four touches $\gamma$, we have $z\not\in \{v_2, v_3,x_2\}$.
If $z=v_1$, then the cycle $K=v_1v_{11}v_8v_7x_7$ separates $z_1$ from $v_2$
\rt{by (I8)}, and thus $x_2\in V(K)$.
This is a contradiction, since then a cycle of length at most four touches $\gamma$.  Therefore, $z$ is
distinct from and non-adjacent to $v_2$.  Furthermore, $z$ is not adjacent to $x_2$, as otherwise
$x_2zx_7\rt{z_1}$ touches $\gamma$.~\qed
\bigskip

Let $G$ be a graph in a surface $\Sigma$ with rings $\cal R$,
let $\gamma$ be a good configuration that weakly appears in $G$,
let $G'$ be the $\gamma$-reduction of $G$, and let $C'$ be a cycle in $G'$.
If $C$ is a cycle in $G$ such that either $C=C'$ or $C'$ is obtained from 
$C$ by replacing a squashed edge by one of the corresponding edges of $G$,
then we say that $C$ is a \emph{lift} of $C'$.

\begin{lemma}
\mylabel{lem:4cycles}
Let $G$ be a graph in a surface $\Sigma$ with rings $\cal R$ satisfying {\rm (I0)}, {\rm (I3), }{\rm (I8)} and {\rm (I9)},
let $\gamma$ be a good configuration that strongly appears in $G$,
and let $G'$ be the $\gamma$-reduction of $G$ with respect to a 
$3$-coloring $\phi$ of $\cal R$.
If $C'$ is a cycle in $G'$ of length at
most four, then either a lift of $C'$ is a cycle in $G$, 
or $C'$ is noncontractible
and there exists a noncontractible cycle $C$ in $G$ such that $C$ touches $\gamma$ and $|C|-|C'|\le3$.
Furthermore, all ring vertices of $C'$ belong to $C$; and if $C'$ is a triangle disjoint from the rings
and its vertices have distinct pairwise non-adjacent neighbors in a ring $R$ of length $6$, then
$G$ contains edges $cr$ and $c'r'$ with $c,c'\in V(C)\setminus V(R)$ and $r,r'\in V(R)\setminus V(C)$
such that $r$ and $r'$ are non-adjacent.
\end{lemma}

\proof 
Suppose that $C'\subseteq G'$ is a cycle of length $3$ or $4$ such that no lift of $C'$ is a cycle in $G$.
Let us discuss the possible configurations $\gamma$:
\begin{itemize}
\item $\gamma$ is isomorphic to one of $\R1$, $\R2$, $\R6$, $\R6.1$, $\R7.1$, $\R7.2$, $\R7.3$ or $\R7.4$, or
to $\R4$ and both $x_4$ and $x_5$ belong to ${\cal R}$ and $\phi(x_4)=\phi(x_5)$.  We are adding an edge $e$ between
vertices $x,x'\in {\cal A}_{\gamma}$ along the replacement path $P\subset G$ of length at most $4$. 
Note that $e\in E(C')$.  Let $C\subseteq G$ be the cycle obtained from $C'$ by replacing $e$ with $P$.
Clearly, $|C|\le |C'|+3\le 7$.  Let us remark that $C$ is indeed a cycle (i.e., if $\gamma$ is $\R4$, then $v_2\not\in V(C')$),
since no non-ring cycle of length at most four touches $\gamma$
by the definition of weak appearance. 
Note that $P$ is not a part of a boundary of a face in
any of the configurations; thus $C$ does not bound a face in $G$.
By (I9), $C$ is not contractible; hence $C'$ is not contractible, either.

\item $\gamma$ is $\R3$:  Let $w$ be the vertex of $G'$ obtained by identifying $v_1$ with $v_3$ and $v_5$.
Note that $w\in V(C')$ and consider the edges $e_1, e_2\in E(C')$ incident with $w$.
Unless $C'$ corresponds to a cycle of length $|V(C')|$ in $G$, 
$e_1$ and $e_2$ are incident with distinct
vertices $a,b\in {\cal I}_{\gamma}$, and the cycle $C$ obtained from $C'$ by adding the replacement path $avb$
between $a$ and $b$ has length at most $|C'|+2\le 6$.  Note that $C'$ and $C$ have the same homotopy.  Suppose that they are contractible.
By (I9) that implies that $C'$ bounds a face $h$ and $v$ has degree two.  By (I0), $v$ belongs to ${\cal R}$;
however, this is not possible, since at least one of $a$ and $b$ is an internal vertex.
This is a contradiction.

\item $\gamma$ is $\R4$ and at least one of $x_4$ and $x_5$ is internal:
Let $w$ be the vertex obtained by identifying $x_4$ and $x_5$.  
If $x_1x_3$ is not an edge of $C'$, then \rt{(since $C'$ is not a cycle of $G$)}
the cycle $C$ obtained from $C'$ by replacing $w$ by the path
$x_4v_4v_5x_5$ satisfies $6\le |C|\le 7$ and
does not bound a face; thus neither $C$ nor $C'$ is contractible.
Let us assume that $x_1x_3\in E(C')$.  Similarly, we deal with the case that $w\not\in V(C')$ or that both edges incident with $w$ in
$C'$ correspond to edges incident to one of $x_4$ and $x_5$.

Suppose now that the neighbors of $w$ in $C'$ are adjacent to $x_4$ and $x_5$.
Since no non-ring cycle of length at most four touches $\gamma$
by the definition of weak appearance, 
we have $x_1x_5,x_3x_4\not\in E(G)$;
thus by symmetry we may assume that $x_1x_4\in E(C')$
and $x_3$ and $x_5$ are joined by a path $P$ of length at most two in $C'$.
By (I8), the $5$-cycle $K=x_1v_1v_5v_4x_4$ separates $x_3$ from $x_5$; 
thus $P$ is not disjoint from $K$.  
However, then a cycle of length at most four touches $\gamma$.

\item $\gamma$ is $\R4$, neither $x_4$ nor $x_5$ is internal and $\phi(x_4)\neq \phi(x_5)$:
Let $w$ be the vertex created by identifying $v_2$ and $x_5$.
The claim of the lemma follows by considering the non-facial cycle $C$ obtained from 
$C'$ by replacing $w$ with $v_2v_1v_5x_5$.

\item $\gamma$ is $\R5$: Let $w$ be the vertex obtained by identifying $v_4$ and $x_6$.  Let $C$ be the cycle
obtained from $C'$ by replacing $v_2x_8$ by $v_2v_1v_8x_8$ or $w$ by $v_4v_5v_6x_6$
or both.  If we performed at most one replacement, then $|C|\le|C'|+3$
and the claim follows from (I9).

Otherwise, $v_2x_8\in E(C')$ and $w\in V(C')$, and since no \rt{non-ring}
cycle of length at most four touches $\gamma$,
there exist paths $P_1$ between $v_2$ and $x_6$ and $P_2$ between $v_4$ and $x_8$ of total length at most three.
Let $K_1$ be the cycle consisting of $v_2v_3v_7v_6x_6$ and $P_1$
and $K_2$ the cycle consisting of $v_4v_3v_7v_8x_8$ and $P_2$, and by symmetry assume that $|K_1|=5$
and $|K_2|\le 6$.  By (I8) the cycle $K_1$ separates $v_4$ from $v_8$; 
thus $P_2$ intersects $K_1$.
However, that contradicts the fact that no non-ring
cycle of length at most four touches $\gamma$.

\item $\gamma$ is $\R7$: Let $w$ be the vertex obtained by identifying $x_6$ and $x_7$.  Let $C_1$ be the cycle obtained
from $C'$ by replacing $x_1x_3$ by $x_1v_1v_2v_3x_3$ or $w$ by $x_6v_6v_7x_7$ or both.
If we performed only one replacement, then $|C_1|=|C'|+3$ and the claim of the lemma follows from (I9), with $C=C_1$.

Otherwise, let $C_2$ be the closed walk obtained from $C_1$ by replacing $x_6v_6v_7x_7$ by $x_6zx_7$;
we have $|C_2|=|C'|+5\le 9$.  Since $\gamma$ appears, observe that all vertices of $C'$ are internal
and at most one of them has a neighbor in a ring.
Note that $C_2$ is a cycle, since otherwise a non-ring 
cycle of length at most four touching~$\gamma$
is a subgraph of $C_2$.  Suppose now that $C_2$ consists of $x_1v_1v_2v_3x_3$, a path $P_1$ from $x_3$ to $x_7$,
the path $x_7zx_6$ and a path $P_2$ from $x_6$ to $x_1$, where the total length of $P_1$ and $P_2$ is at most three.
Let $K_1$ be the cycle consisting of $P_1$ and $x_3v_3v_{12}v_6v_7x_7$
and $K_2$ the cycle consisting of $P_2$ and $x_1v_1v_{11}v_{12}v_6x_6$.
Note that $\min(|K_1|,|K_2|)\le 6$, and by (I8), the shorter of the two cycles is separating.
It follows that $K_1$ and $K_2$ intersect in a vertex distinct
from $v_{12}$ and $v_6$.  This is a contradiction, since the vertices of $C_2$ are
mutually distinct and none of them is equal to $v_7,v_{11}\not\in V(G')$.

Therefore, $C_2$ consists of $x_1v_1v_2v_3x_3$, a path $Q_1$ of length $l_1\ge 1$ from $x_3$ to $x_6$,
the path $x_6zx_7$ and a path $Q_2$ of length $l_2$ from $x_7$ to $x_1$, where $l_1+l_2\le 3$.
Let $L_1$ be the cycle consisting of $Q_1$
and $x_3v_3v_{12}v_6x_6$ and $L_2$ the cycle consisting of $Q_2$ and $x_1v_1v_{11}v_8v_7x_7$.
Note that neither $L_1$ nor $L_2$ bounds a face, $|L_1|=4+l_1\le 7$ and
$|L_2|=5+l_2\le 7$, thus by (I9) neither $L_1$ nor $L_2$ is contractible.
Furthermore, $|L_1|+|L_2|\le 9+l_1+l_2\le 12$, thus there exists a cycle $C\in \{L_1,L_2\}$ of
length at most $6\le |C'|+3$ touching~$\gamma$.

Let us now show that the cycle $C'$ is not contractible.
Assume for a contradiction that $C'$, and hence also $C_2$, is contractible.
Let $\Delta\subseteq \Sigma$ be an open disk bounded by $C_2$.  Note that $\Delta$ does not consist of a single face,
since at least one edge incident with $v_1$ or
$v_2$ lies inside $\Delta$.  By (I9), $\Delta$ consists of two or three faces, and in the latter case,
$|C_2|=9$ and three vertices of $C_2$ have a common neighbor.  

It follows that $v_{11},v_{12}\not\in\Delta$, and thus the edge joining $v_2$ with its neighbor $x_2\not\in\{v_1,v_3\}$
lies in $\Delta$.  Since $\gamma$ appears strongly in $G$, we have that $x_2\neq z$ and that $z$ is an internal vertex.
We conclude that $\deg(z)=3$ and $z$ has a neighbor inside $\Delta$ distinct from $x_6$ and $x_7$. 
By (I3) and (I9), this neighbor is equal to $x_2$.  However, this contradicts the assumption that $\gamma$ appears strongly in $G$.
\end{itemize}

If $\gamma$ is $\R7$ and $C$ is one of the cycles $L_1$ and $L_2$, then since $\gamma$ appears in $G$,
the vertices $x_1$, $x_3$, $x_6$, $x_7$ and all their neighbors in $G$ are
internal.  Consequently, $x_1$, $x_3$ and all their neighbors are internal in $G'$.  It follows that $C'$ contains no ring vertex,
and that at most two distinct ring vertices have a neighbor in $C'$, hence the last claim of the lemma holds trivially.

Otherwise, $C$ is obtained from $C'$ by replacing a new edge by a path in $G$, or by adding a replacement path between vertices
of ${\cal I}_{\gamma}$, or both.  Therefore, every ring vertex of $C'$ also belongs to $C$.  Suppose that $C'$ is a triangle
whose vertices $c_1$, $c_2$ and $c_3$ are internal, that $R=r_1r_2r_3r_4r_5r_6$ is a ring and that $c_1r_1$, $c_2r_3$ and $c_3r_5$ are edges
of $G'$.  If, say, $r_1$ has no neighbor in $C$, then either $r_1c_1$ is a new edge, or one of $r_1$ and $c_1$ is the new vertex created by the identification of
the vertices of ${\cal I}_{\gamma}$.  Since $C$ is not a lift of $C'$, in the former case $C'$ contains a new vertex that is replaced by a path in $C$,
and in the latter case $C'$ contains a new edge. Therefore, $c_2r_3$ and $c_3r_5$ are edges of $G$.

\qed

\section{Contributions of faces}\label{sec-winners}

Let $G$ be a graph in a surface $\Sigma$ with rings $\cal R$ satisfying (I3).
Let $\gamma$ be a good configuration that strongly appears in $G$,
let $G'$ be the $\gamma$-reduction of $G$,
and let $G''$ be a subgraph of $G'$ that includes all the rings and satisfies (I0).

Let $f''$ be a face of $G''$, and let $H$ be the subgraph
of $G''$ that forms the boundary of $f''$.
We wish to define a subgraph $J_{f''}$ of $G$ that will correspond to $H$,
and a union of faces of $J_{f''}$ that will correspond to $f''$.

Let us recall that during the construction of the graph $G'$, 
parallel edges may
have been removed (e.g., if $\gamma$ is $\R5$ and $v_4$ and $x_6$ have a common neighbor), but we have retained the correspondence
of each non-squashed edge $e$ of $G'$ to a unique edge of $G$ (which also determined the placement of $e$ in the embedding of $G'$).
We now define the edge-set of $J_{f''}$, by replacing pieces of the boundary of $f''$ by appropriate replacement
paths.  More precisely, we apply the following construction to each boundary walk $C$
of $f''$.  Let $C$ be $v_1, e_1, v_2, e_2, \ldots, v_m, e_m$
and let $e_{m+1}=e_1$, $v_{m+1}=v_1$, $e_0=e_m$ and $v_0=v_m$.
Suppose that $v_i$ is a vertex and $e\in\{e_{i-1}, e_i\}$ is an edge of $C$ incident with $v_i$, and let $e'\in \{e_{i-1}, e_i\}\setminus \{e\}$ be the other
edge of $C$ incident with $v_i$.  Note that $e_{i-1}\neq e_i$, since otherwise $v_i$ would have degree $1$ in $G''$,
contrary to the assumption that $G''$ includes all rings and satisfies (I0).  We define $\text{orig}(v_i,e)$ as follows.
\begin{itemize}
\item If $v_i$ is not a new vertex, then $\text{orig}(v_i,e)=v_i$.
\item If $v_i$ is a new vertex and $e$ is not a squashed edge, then the edge of $G$ corresponding to $e$ is incident with a unique vertex
$z\in {\cal I}_\gamma$, and we define $\text{orig}(v_i,e)=z$.
\item If $v_i$ is a new vertex and both $e$ and $e'$ are squashed edges, then the inspection of configurations shows that this is only possible
if $\gamma$ is isomorphic to the configuration $\R3$.  In this case, $\text{orig}(v_i,e)$ is defined to be the vertex which is in the depiction
of $\R3$ in Figure~\ref{fig-redu} denoted by $v_5$.
\item Finally, suppose that $v_i$ is a new vertex, $e$ is a squashed edge with the other end $u\neq v_i$, and $e'$ is not a squashed edge.
If $u$ is adjacent to $\text{orig}(v_i,e')$, then $\text{orig}(v_i,e)=\text{orig}(v_i,e')$.  Otherwise, note that $\gamma$ is isomorphic to the configuration $\R3$,
and $\text{orig}(v_i,e)$ is again defined to be the vertex which is in the depiction
of $\R3$ in Figure~\ref{fig-redu} denoted by $v_5$.
\end{itemize}
Define $\text{orig}(e_i)$ as the edge of $G$ joining $\text{orig}(v_i,e_i)$ with $\text{orig}(v_{i+1},e_i)$; note that $\text{orig}(e_i)=e_i$ unless $e_i$ is
a squashed edge.

Now, replace each edge $e_i$ of $C$ by a path $P_i$ defined as follows.
If $e_i$ is a new edge, then $P_i$ is the corresponding replacement path.
Otherwise, let $P_i$ consist of the edge $\text{orig}(e_i)$, and in case that $\text{orig}(v_{i+1},e_i)\neq \text{orig}(v_{i+1},e_{i+1})$
also of the replacement path between the vertices $\text{orig}(v_{i+1},e_i)$ and $\text{orig}(v_{i+1},e_{i+1})$ (the two vertices belong to ${\cal I}_\gamma$).
The newly constructed walk has the same homotopy as $C$.
The graph $J_{f''}$ is defined as the result of applying the above
construction to every boundary walk of $f''$.

\begin{figure}
\begin{center}
\epsfxsize=12cm\epsfbox{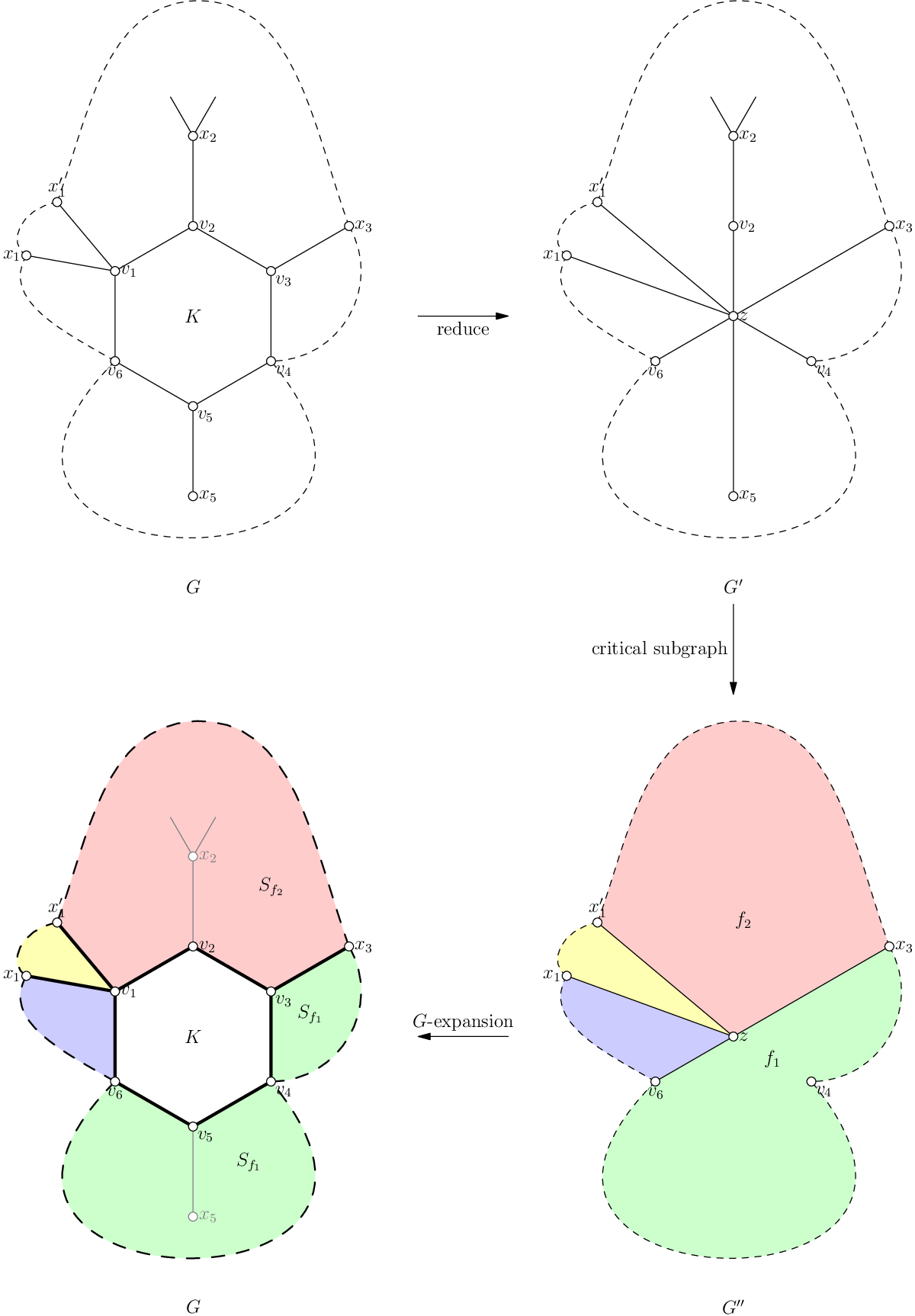}
\end{center}
\caption{The $G$-expansion.}\label{fig-mainstep}
\end{figure}

The construction is illustrated in Figure~\ref{fig-mainstep}, where configuration $\R3$
is reduced.  For example, in the boundary walk $v_6zx_3\ldots v_4\ldots$ of the
face $f_1$ of $G''$, we have $\text{orig}(z,v_6z)=v_5$ and $\text{orig}(z,zx_3)=v_3$,
and thus the edge $v_6z$ is replaced by $\text{orig}(v_6z)=v_6v_5$ and the replacement
path $v_5v_4v_3$.

As this example also illustrates, a face $f''$ of $G''$,
it may correspond to several faces of $J_{f''}$, in case that an internal vertex of a
replacement path belongs to the subgraph $G''$; let the set of these faces
of $J_{f''}$ be denoted by $S_{f''}$.

The {\em elasticity} of $f''$ is $\el(f'')=\left(\sum_{f\in S_{f''}}|f|\right)-|f''|$.
In the example, $\el(f_1)=\el(f_2)=2$ and $\el(f_3)=\el(f_4)=0$.
Note that $f''$ can have non-zero elasticity only if $J_{f''}$ contains
at least one replacement path and each replacement path contributes to elasticities
of at most two faces.  Furthermore, if the path contributes to the elasticity of a face $f''$ twice,
then the corresponding new vertex or both vertices incident with the corresponding new edge
appear at least twice in the boundary walks of $f''$.  
This cannot happen if $f''$ is semi-closed $2$-cell and all vertex-like rings of $G''$ are also vertex-like in $G$, since a new vertex cannot be the main vertex of a vertex-like ring
by the first condition in the definition of appearance and the analogous later condition for $\R4$, and the new edge
cannot join the main vertices of distinct vertex-like rings by the first condition in the definition of strong appearance.
Using these observations and the inspection of the configurations, we obtain the following.

\begin{lemma}
\mylabel{lem:elasticity}
Let $G$, $\gamma$, $G'$, and $G''$ be as above.  Then $G''$ has at most three faces with
non-zero elasticity, and the sum of the elasticities of the faces of $G''$ is at most $10$.
If a face $f''$ of $G''$ is closed $2$-cell or omnipresent, then $\el(f'')\le5$, and if the inequality is strict, then $\el(f'')\le3$.
Furthermore, if all vertex-like rings of $G''$ are also vertex-like in $G$, then the previous statement holds also
for semi-closed $2$-cell faces.
\end{lemma}

Let $G$ be a graph in a surface $\Sigma$ with rings $\cal R$, let
$J$ be a subgraph of $G$, and let $S$ be a subset of the set of
faces of $J\cup\bigcup{\cal R}$ such that $J$ is equal to the union of the boundaries of the faces in $S$.
We define $G[S]$ to be the subgraph of $G$ consisting of $J$
and all the vertices and edges drawn inside the faces of $S$.
Let $C_1,C_2,\ldots,C_k$ be the boundary walks of the faces in $S$.
We would like to view $G[S]$ as a graph with rings $C_1$, \ldots, $C_k$.
However, the walks $C_1$, \ldots, $C_k$ do not necessarily have to be disjoint, and they do not have to be
cycles.
To overcome this difficulty, we proceed as follows:
Suppose that $S=\{f_1,\ldots, f_m\}$.  For $1\le i\le m$, let $\Sigma_i$ be
a surface with boundary $B_i$ such that $\Sigma_i\setminus B_i$ is
homeomorphic to $f_i$. Let $\theta_i:\Sigma_i\setminus B_i\to f_i$ be a homeomorphism 
that extends to a continuous mapping 
$\theta_i:\Sigma_i\to\overline{f_i}$, where $\overline{f_i}$ denotes the closure of $f_i$.
Let $G_i$ be the inverse image of $G\cap \overline{f_i}$ under $\theta_i$.
Then $G_i$ is a graph normally embedded in $\Sigma_i$.
We say that the set of embedded graphs $\{G_i:1\le i\le m\}$ 
is a {\em $G$-expansion of $S$.}
Note that there is a one-to-one correspondence between the boundary walks of the faces of $S$ and the rings of the graphs
in the $G$-expansion of $S$; however, each vertex of $J$ may be split to several copies.

We define the \emph{$G$-expansion of $f''$} to be the $G$-expansion of $S_{f''}$.  The following lemma is straightforward.

\begin{lemma}
\mylabel{lem:facecover}
Let $G$, $\gamma$, $G'$, and $G''$ be as above, and let $f$ be a face of $G$.
Then either there exists a unique face $f''$ of $G''$ such that
$f$ corresponds to a face of a member of 
the $G$-expansion of $f''$ or $\gamma$ is isomorphic to $\R3$ and
$f$ is the $6$-face of ${\cal F}_\gamma$.
\end{lemma}

Let us now give an informal summary of what we are trying to achieve in this section.
We assign weights to the faces of embedded graphs according to the function $s$
as described in Section~\ref{sec-disch},
and we aim to show that the sum of the weights of the faces of $G$ is
bounded by the sum of the weights of the faces of $G''$.  To do so, we would
like to claim that the sum $w$ of the weights of the faces 
of members of the $G$-expansion of $f''$ is bounded by
the weight $w''$ of $f''$.  
In Theorem~\ref{thm:diskgirth5}, we will show that this
claim holds, provided that the elasticity of $f''$ is small and 
the $G$-expansion of $f''$ is not a singleton set consisting of
one of a few exceptional graphs.  Here, we assign a {\em contribution} $c(f'')$ to each face $f''$ of $G''$
according to the criteria that we later prove to ensure that $w\le w''-c(f'')$.  Furthermore, we argue that
the sum of the contributions of all faces is non-negative.

Let us now proceed more formally.  
Let $G$ be a graph in a surface $\Sigma$ with rings $\cal R$.
If $\Sigma$ is a disk and ${\cal R}=\{R\}$, then we say that
$G$ is a \emph{plane graph with one ring $R$}.
We say that a plane graph $G$ with one ring $R$ of length $l\ge5$
is {\em exceptional} if it satisfies one of the conditions below (see Figure~\ref{fig-except}):

\begin{figure}
\epsfbox{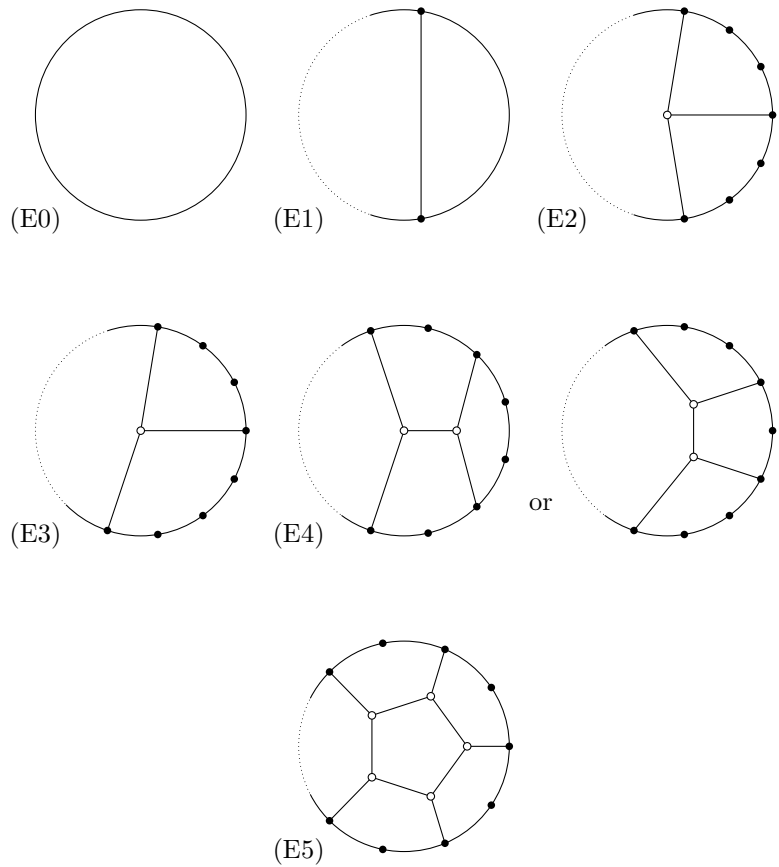}
\caption{Exceptional graphs.}
\label{fig-except}
\end{figure}

\begin{itemize}
\item[(E0)] $G=R$,
\item[(E1)] $l\ge 8$ and $|E(G)|-|E(R)|=1$,
\item[(E2)] $l\ge 9$, $V(G)-V(R)$ has exactly one vertex of degree three,
and the faces of $G$ have lengths $5,5,l-4$,
\item[(E3)] $l\ge 11$, $V(G)-V(R)$ has exactly one vertex of degree three,
and the faces of $G$ have lengths $5,6,l-5$,
\item[(E4)]  $l\ge 10$, $V(G)-V(R)$ consists of two adjacent degree three vertices,
and the faces of $G$ have lengths $5,5,5,l-5$,
\item[(E5)] $l\ge 10$, $V(G)-V(R)$ consists of five degree three vertices forming
a facial cycle of length five,
and the faces of $G$ have lengths $5,5,5,5,5,l-5$.
\end{itemize}
\noindent We say that $G$ is {\em very exceptional} if it satisfies
(E0), (E1), (E2) or (E3).

Let us now show the following lemma, which we use to analyze omnipresent faces.

\begin{lemma}\label{lemma-rednt}
Let $G$ be a graph in a surface $\Sigma$ with rings $\cal R$,
let $G$ be well-behaved, let it satisfy (I0), (I4) and (I8),
let $\gamma$ be a good configuration that strongly appears in $G$, 
let $G'$ be the $\gamma$-reduction of $G$,
let $G''$ be a subgraph of $G'$ that includes all the rings and
satisfies (I0),
and let $H$ be a component of $G''$ that contains a new edge or a new vertex.  
Assume that either $\Sigma$ is a disk and $|{\cal R}|=1$ and every face of 
$G''$ is closed $2$-cell, 
in which case we let $R$ be the unique member of $\cal R$,
or that $G''$ has an omnipresent face, 
in which case we let $R$ be the boundary walk of the omnipresent face 
that is a subgraph of $H$.  
In either case $H$ can be regarded as a plane graph with one ring $R$.
Then $H$ is not very exceptional.
Furthermore, if $\gamma$ is isomorphic to one of $\R6$, $\R6.1$, $\R7$, $\R7.1$, $\R7.2$, $\R7.3$ or $\R7.4$, then $H$ is not exceptional.
\end{lemma}

\begin{proof}
Since $H$ contains a new vertex or a new edge (which are not contained in the boundary), it does not satisfy (E0).
If $\gamma$ is isomorphic to one of $\R7$, $\R7.1$, $\R7.2$, $\R7.3$ or $\R7.4$, then all vertices in ${\cal A}_\gamma\cup{\cal I}_\gamma$ and all their neighbors
are internal by the definition of appearance, 
and thus each new edge or new vertex is at distance at least two from $R$.
It follows that $H$ cannot be exceptional.  Similarly, we exclude the case that $\gamma$ is isomorphic to $\R6$ or $\R6.1$.
Thus, assume that $\gamma$ is one of $\R1$, \ldots, $\R5$.

Suppose that $H$ contains a new edge $xy$.  
Note that since $\gamma$ does not touch a non-ring cycle of length at most 
four by the definition of weak appearance,
neither $x$ nor $y$ is a new vertex.  Since $\gamma$ appears strongly in $G$, 
we may assume that $x$ is an internal vertex;
thus $H$ does not satisfy (E1).  Suppose that $H$ satisfies (E2) or (E3).  
Then, in $H$ the vertex $x$ has three neighbors in $R$.  
On the other hand, $x$ has at most one neighbor in $R$ in $G$, by (I4).  We conclude that $x$ is adjacent to a new vertex in $G''$
that belongs to $R$.  It follows that $\gamma$ is $\R4$ or $\R5$, and in the former case at least one of $x_4$ and $x_5$ is internal.
Let ${\cal I}={\cal I}_\gamma$ if $\gamma$ is $\R5$ and ${\cal I}=\{x_4,x_5\}$ if $\gamma$ is $\R4$.  Note that there exists a vertex
in ${\cal I}$ belonging to $R$, and another vertex of ${\cal I}$ is adjacent to $x$ in $G$.
If $\gamma$ is isomorphic to $\R4$, then by symmetry we may assume that $x_1$ is adjacent to $x_4$ and
$x_3$ and $x_5$ belong to $R$.  However, by (I8), the cycle $x_1v_1v_5v_4x_4$ consisting of internal vertices
separates $x_3$ from $x_5$, which is a contradiction.
If $\gamma$ is isomorphic to $\R5$, then by the conditions of appearance, $v_4$ is an internal vertex; hence $x_6$ belongs to $R$.
Since $v_2$ and $v_4$ are not adjacent, we conclude that $v_4$ is adjacent to $x_8$ and that $v_2$ belongs to $R$.
However, this again contradicts (I8).

Therefore, we may assume that $H$ contains a new vertex, but not a new edge.  Suppose first that $\gamma$ is not
isomorphic to $\R4$.  If $H$ satisfied (E1), then by (I4) there would exist vertices $x\in {\cal I}_\gamma\cap V(R)$
and $y\in {\cal I}_\gamma\setminus V(R)$ and a neighbor $z$ of $y$ in $R$, where $z$ is not adjacent to $x$.  However,
this contradicts the assumption that $\gamma$ appears strongly in $G$.
If $H$ satisfies (E2) or (E3), then by (I4) we have $|{\cal I}_\gamma|=3$ (thus $\gamma$ is $\R3$), all elements of
${\cal I}_\gamma$ are internal and each of them has exactly one neighbor in $R$.  This is excluded, since $\gamma$
appears in $G$.

Finally consider the case that $\gamma$ is $\R4$ and $H$ does not contain a new edge.  By (I4), $H$ does not satisfy
(E2) or (E3); thus suppose that $H$ satisfies (E1).  If $x_4$ is an internal vertex, this implies that $x_5\in V(R)$
and $x_4$ has a neighbor $w$ in $R$ distinct from $z$.  By (I4), $z$ is an internal vertex.  Since $G$ is well-behaved,
the path $x_5zx_4w$ forms a part of a boundary of a $5$-face; thus $z$ has degree two, contrary to (I0).
The case that $x_5$ is internal is symmetric; thus assume that both $x_4$ and $x_5$ belong to $R$.
Then $v_2$ is an internal vertex of degree at least four 
by the definition of weak appearance and has a neighbor $w\in V(R)$.  
However, since $G$ is
well-behaved, the subpaths $v_3v_2w$ and $v_1v_2w$ of the paths $x_4v_4v_3v_2w$ and $x_5v_5v_1v_2w$ form parts
of boundaries of faces, implying on the contrary that $v_2$ has degree three.
\end{proof}

Let $s:\{5,6,\ldots\}\to{\mathbb R}^+$ be an increasing function, to be 
specified later, such that

{
\myclaim{S4} {$14s(5)\le s(6)$, $135s(5)\le s(7)$, $4s(6)\le s(7)$, 
$3s(7)\le s(8)$, $2s(8)\le s(7)+s(9)$ and $s(l)=l-8$ for $l\ge 9$.}
}

\noindent
It follows that the function $s$ satisfies
\ppar
{
\myclaim{S5} {$s(x+a)-s(x)\le s(y+a)-s(y)$ for all integers $y\ge x\ge5$
and $a\ge0$.}
}

\noindent
We will refer to condition (S5) as \emph{convexity}.

Let $G$ be a graph in a surface $\Sigma$ with rings $\cal R$,
let $\gamma$ be a good configuration that strongly appears in $G$,
let $G'$ be the $\gamma$-reduction of $G$, and let $G''$ be a subgraph of
$G'$ that includes all the rings and satisfies (I0).
For every face $f''$ of $G''$ we define its \emph{contribution}
$c(f'')$ as follows.

Let $f''$ be a semi-closed $2$-cell face of $G''$, 
and let $G_{f''}$ be a member of the $G$-expansion of $f''$. 
Then the contribution of $f''$ is defined 
according to the following rules:
\begin{itemize}
\item If $|S_{f''}|=1$ and
$G_{f''}$ satisfies (E0), then $c(f'')=-\infty$ if $f''$ has non-zero elasticity and $c(f'')=0$ otherwise.
\item If  $|S_{f''}|=1$ and
$G_{f''}$ satisfies (E1), then $c(f'')=-\infty$ if $\el(f'')=5$ and $c(f'')=s(8-\el(f''))-2s(5)$ otherwise.
\item If  $|S_{f''}|=1$ and
$G_{f''}$ satisfies (E2), then $c(f'')=-\infty$ if $\el(f'')=5$ and $c(f'')=s(9-\el(f''))-3s(5)$ otherwise.
\item If  $|S_{f''}|=1$ and
$G_{f''}$ satisfies (E3), then $c(f'')=s(11-\el(f''))-2s(6)-s(5)$.
\item If  $|S_{f''}|=1$ and
$G_{f''}$ satisfies (E4) or (E5), or if $|S_{f''}|=2$ and $G_{f''}$ consists of two cycles such that one of them has length $5$,
then $c(f'')=s(10-\el(f''))-6s(5)$.
\item If $|S_{f''}|=1$ and $G_{f''}$ is not exceptional, and 
\begin{itemize}
\item $G_{f''}$ contains a path $P=v_1v_2v_3v_4$ such that $v_1,v_4\in V(J_{f''})$, $v_2,v_3\not\in V(J_{f''})$ and both of the open disks
bounded by $P$ and paths in $J_{f''}$ contain at least two vertices of $G$, then $c(f'')=s(7)$.
\item Otherwise, $c(f'')=s(11-\el(f''))-s(6)+5s(5)$.
\end{itemize}
\item If $|S_{f''}|=2$ and $G_{f''}$ does not consist of two cycles such that one of them has length $5$, or if $|S_{f''}|\ge 3$, then
$c(f'')=s(12-\el(f''))-2s(6)$.
\end{itemize}

Suppose now that $f''$ is an omnipresent face of $G''$.  Let $G''_1$, $G''_2$, \ldots, $G''_k$
be the components of $G''$ such that $G''_i$ contains the ring $R_i\in{\cal R}$.  If there exist $i\neq j$
such that $G''_i\neq R_i$ and $G''_j\neq R_j$, then $c(f'')=1$.  
Otherwise, we may assume that $G''_i=R_i$ for $i\ge 2$.
If $G''_1$ satisfies (E0), (E1), (E2) or (E3), then $c(f'')=-\infty$.  If $G''_1$ satisfies (E4) or (E5),
then $c(f'')=5-\el(f'')-5s(5)$, otherwise $c(f'')=5-\el(f'')+5s(5)$.

This completes the definition of contribution of faces of $G''$.
We define the contribution of $G''$ as
$c(G'')=-\delta+\sum_{f''\in F(G'')} c(f'')$, where $\delta$ is $s(6)$ if $\gamma$ is isomorphic to $\R3$
and $0$ otherwise.

\begin{lemma}
\mylabel{lem:winners}
Let $G$ be a well-behaved graph in a surface $\Sigma$ with rings $\cal R$
satisfying (I0)--(I4) and (I8),
let $\gamma$ be a good configuration strongly appearing in $G$, and
let $G'$ be the $\gamma$-reduction of $G$.
Suppose that $G''$ is a subgraph of $G'$ that includes $\cal R$,
every vertex-like ring of $G''$ is also vertex-like in $G'$,
$G''$ satisfies (I0) and (I6), and
$G''$ contains a new vertex or a new edge.
\begin{enumerate}
\item If each face of $G''$ is semi-closed $2$-cell or omnipresent, then $c(G'')\ge 0$.
\item If each face of $G''$ is semi-closed $2$-cell, then $G''$ has a face of
length at least six.
\item If each face of $G''$ is semi-closed $2$-cell, $\Sigma$ is a disk, and $|{\cal R}|=1$, then $c(G'')\ge 10s(5)$.
\end{enumerate}
\end{lemma}

\proof
We first prove the first and third statements of the lemma.
The proof will show that there is a face of positive contribution,
and we will utilize that face in the proof of the second statement of the lemma.
Let us note that $G''$ satisfies the assumptions of Lemma~\ref{lemma-rednt}, and thus
the contribution of the omnipresent face (if $G''$ has such a face) is not $-\infty$.

We may assume that there exists a face of non-zero elasticity,
for otherwise all faces have non-negative contribution, and
all vertices of $\dom(d_\gamma)$ are included in the interior of
a single face $f''$ of $G''$; clearly, $G_{f''}$ is not exceptional, and thus
this face satisfies $c(f'')\ge s(7)$.

Let us argue that if a face $f''$ that is semi-closed $2$-cell
has non-zero elasticity, then either $S_{f''}$ has at least
two faces or the interior of the unique
face of $S_{f''}$ contains an edge of $G$.  Indeed, most replacement
paths are incident with edges on both of its sides; thus if such a replacement
path is used in $J_{f''}$, then at least one such edge lies in $S_{f''}$.
The exceptions are the replacement paths in $\R3$, $\R4$ and the replacement
path between the vertices of ${\cal I}_\gamma$ in $\R7$.  
In these configurations,
the middle vertex $v$ of the replacement path could also lie on the boundary walk of $f''$, in which case all
the edges incident with $v$ could belong to $J_{f''}$ or lie outside of $S_{f''}$.
However, then $S_{f''}$ has at least two faces.  We conclude that
if $c(f'')=-\infty$, then $\el(f'')=5$ and two replacement paths are used in the construction of $J_{f''}$.

Let us first consider the case when every face of $G''$ has elasticity 
at most three. Then
the contribution of each face is greater or equal to $-s(5)$.
Furthermore, there are at most two faces of elasticity three, 
at most one of them has contribution $-s(5)$, and
every other face has non-negative contribution.
If $G''$ has an omnipresent face, then $c(G'')\ge 2-6s(5)-s(6)$
or $c(G'')\ge 1-s(6)$, and 
hence we may assume that all faces of $G''$ are semi-closed $2$-cell.
Observe that $G$ has a face $f''$ such that at least one vertex of
$\dom(d_\gamma)$ is contained inside a face of $S_{f''}$.
For this face, we have $c(f'')\ge s(6)-3s(5)$.  
Furthermore, if $\gamma$ is $\R3$, then
the elasticity of $f''$ is at most two; thus $c(f'')\ge s(7)-3s(5)$, and all other faces of $G''$ have non-negative contribution.  
Therefore, $c(G'')\ge \min(s(6)-4s(5),s(7)-s(6)-3s(5))\ge 10s(5)$,
because $14s(5)\le s(6)$.
This completes the case when every face of $G''$ has elasticity at most three.

Thus we may assume that $G''$ has a face $f_2$ of elasticity five.
It follows that $\gamma$ is $\R4$, $\R5$, or $\R7$,
and $G''$ contains a new edge incident with two faces of non-zero elasticity, say $f_1$ and $f_2$,
where $f_2$ contains all vertices of $\dom(d_\gamma)$.  Furthermore, $G''$ contains
a new vertex $w$ incident with $f_2$ and possibly another face $f_3$ of non-zero elasticity.

Then the elasticity of $f_2$ is $5$, and by  inspection of the configurations, we conclude that
$c(f_2)\ge -5s(5)$.  Furthermore, if $\gamma$ is isomorphic to $\R7$, then $c(f_2)=s(7)$ if $f_2$ is semi-closed $2$-cell,
and by Lemma~\ref{lemma-rednt}, we have $c(f_2)\ge5s(5)$ if $f_2$ is omnipresent.

Assume now that either $f_2$ is the only face of $G''$ with non-zero elasticity that is incident
with $w$, or that $f_1\neq f_3$.  Consider a face $f\in\{f_1,f_3\}$ with non-zero elasticity.
Since $\el(f)\le 3$, we have $c(f)\ge -s(5)$.  
Furthermore, if $f$ is omnipresent, then by Lemma~\ref{lemma-rednt}, 
we have $c(f)\ge 2-5s(5)$ or $c(f)\ge 1$, and hence
 $c(G'') > 10s(5)$; thus we may  
assume that each such face $f$ is semi-closed $2$-cell.
If $\gamma$ is $\R5$, then $\el(f_1)=2$ and $c(f_1)\ge s(6)-2s(5)$.
Similarly if $\gamma$ is $\R4$, then by 
the definition of appearance $v_2$ has degree at least $4$ in $G_{f_1}$;
hence $c(f_1)\ge s(7)-6s(5)$.  In both cases we get $c(G'')\ge 10s(5)$.  
Thus we may assume that $\gamma$ is $\R7$.
If $\Sigma$ is a disk and $|{\cal R}|=1$,
then $f_2$ is not omnipresent, and hence $c(G'')\ge s(7)-2s(5)\ge 10s(5)$,
because (using the numbering of the vertices as in $\R7$) the path
$v_3v_{12}v_6x_6$ shows that the contribution of $f_2$ is $s(7)$.  
Otherwise, $c(G'')\ge 3s(5)$,
because $c(f_2)\ge 5-\el(f_2)+5s(5)=5s(5)$, using Lemma~\ref{lemma-rednt}.

Therefore, we may assume that $f_1=f_3$ and $f_1$ has elasticity $5$.  If $\Sigma$ were a disk and $|{\cal R}|=1$,
or if $f_1$ or $f_2$ were omnipresent, then $w$ together with a vertex of the new edge would form a $2$-cut
in $G''$, contradicting (I6).  We conclude that both $f_1$ and $f_2$ are semi-closed $2$-cell and that either $\Sigma$ is not
a disk or $|{\cal R}|\neq 1$; hence, it suffices to show that $c(G'')$ is non-negative.

Suppose that $\gamma$ is $\R4$.
Since $\gamma$ weakly appears in $G$, we have that no cycle of length at most $4$ touches $\gamma$, and thus $z\neq v_2$.
The fact that $v_2$ has degree at least four in $G_{f_1}$ implies that $c(f_1)\ge 5s(5)$, unless $G_{f_1}$ consists of a $5$-cycle and a $|f_1|$-cycle.
However, in that case $z$ would be a vertex of degree two, contradicting (I0). It follows that, $c(G'')\ge 0$.

Assume next that $\gamma$ is $\R5$.  By (I1) and (I2) we have that $G_{f_1}$
is not an exceptional graph (considering the cycle formed by the path $v_1v_8v_7v_6v_5$
together with a path in $G_{f_1}$), thus again $c(f_1)\ge 5s(5)$ and $c(G'')\ge 0$.

Finally, let $\gamma$ be $\R7$.  If $|S_{f_1}|\ge 2$, then $c(f_1)\ge -5s(5)$.
Otherwise, note that $z$ is by (I0) incident with an edge lying
inside the face of $S_{f_1}$.  Since $\gamma$ appears strongly in $G$, we have that $v_2$ is not adjacent to $z$, and $v_2$ and $z$
have no common neighbor distinct from $v_1$, $v_3$, $x_6$ and $x_7$. It follows that $G_{f_1}$ does not satisfy
(E1), (E2) or (E3), and thus $c(f_1)\ge -5s(5)$.  Therefore, $c(G'')\ge s(7)-5s(5)>0$.

Therefore, both inequalities from the statement of the lemma hold.  Furthermore,
note that in all the cases, at least one face $f''$ of $G''$ has positive contribution; and if $f''$ is semi-closed $2$-cell, then $|f''|\ge 6$.~\qed

\section{Plane graphs with one ring}
\mylabel{sec-disk}

Before we turn our attention to plane graphs with one ring, let us show several properties of critical
graphs.
Let us recall that $\cal R$-critical graphs were defined at the
end of Section~\ref{sec:def}.

\begin{lemma}
\mylabel{lem:i012}
Let $G$ be a graph in a surface $\Sigma$ with rings $\cal R$.
If $G$ is $\cal R$-critical, then it satisfies {\rm (I0), (I1) and (I2)}.
\end{lemma}

\proof If $G$ contains an internal vertex $v$ of degree at most two, then let $G'=G-v$.
If $G$ contains a cycle $C$ consisting of internal vertices of degree three that has
even length or two vertices of $C$ have adjacent neighbors, then let $G'=G-V(C)$.
For any precoloring $\psi$ of $\cal R$ that extends to a $3$-coloring $\phi$ of $G'$,
observe that $\phi$ can be extended to a $3$-coloring of $G$.  This contradicts the assumption that $G$ is
$\cal R$-critical.~\qed

\bigskip

By the theorem of Gr\"otzsch, no component of a critical graph is a triangle-free planar graph.  This observation can be
strengthened as follows.

\begin{lemma}
\mylabel{lem:crit3conn}
Let $G$ be a graph in a surface $\Sigma$ with rings $\cal R$.  Suppose that each component of $G$ is a planar graph
containing exactly one of the rings.  If $G$ is ${\cal R}$-critical and contains no non-ring triangle,
then each component of $G$ is $2$-connected and $G$ satisfies (I6).
\end{lemma}

\proof We can consider each component of $G$ separately; thus assume that
$\Sigma$ is the sphere and $G$ has only one ring $R$.  Firstly, observe that $G$ is $2$-connected; otherwise, it contains
proper subgraphs $G_1$ and $G_2$ such that $G=G_1\cup G_2$ and $|V(G_1)\cap V(G_2)|\le 1$.  Since $R$ is $2$-connected,
we can assume that $R\subseteq G_1$. However, $G_2$ is $3$-colorable
by Theorem~\ref{grotzsch}, and since we can permute the colors arbitrarily, any precoloring
of the common vertex of $G_1$ and $G_2$ extends to a $3$-coloring of $G_2$.  It follows that any $3$-coloring of $G_1$ extends to
a $3$-coloring of $G$, contrary to the criticality of $G$.

Suppose now that $G$ has an internal $2$-cut, i.e., there exist proper induced subgraphs $G_1$ and $G_2$ of $G$ such that
$G=G_1\cup G_2$, $V(G_1)\cap V(G_2)=\{u,v\}$ for some vertices $u,v\in V(G)$, and $R\subseteq G_1$.  Since $G$ is $2$-connected and planar,
both $u$ and $v$ are incident with the same face of an embedding of $G_2$ in the plane.  If $u$ and $v$ were adjacent, then we would argue as in the
previous paragraph that every precoloring of $u$ and $v$ by distinct colors extends to a $3$-coloring of $G_2$, contrary to the criticality of $G$.
If $u$ and $v$ are not adjacent, then let $G'_2$ be the graph obtained from $G_2$ by adding vertices $z_1$ and $z_2$
and edges of paths $uz_1v$ and $uz_2v$.  The resulting graph is triangle-free, and by \cite{Gro},
every precoloring of the cycle $uz_1vz_2$ using three colors extends to a $3$-coloring of $G'_2$; hence, every precoloring of $u$ and $v$
extends to a $3$-coloring of $G_2$.  Again, this contradicts the criticality of $G$.~\qed

\bigskip

If $G$ is a plane graph with one ring $R$, then we abbreviate
$\{R\}$-critical to {\em $R$-critical}.  Such graphs are very
important for the study of critical graphs on surfaces, for the following
reason:

\begin{lemma}
\mylabel{lem:diskcritical}
Let $G$ be a graph in a surface $\Sigma$ with rings $\cal R$,
and assume that $G$ is $\cal R$-critical.
Let $C$ be a non-facial cycle in $G$ bounding an open disk 
$\Delta\subseteq \Sigma$, and let $G'$ be the graph
consisting of the vertices and edges of $G$ drawn in the closure of $\Delta$.
Then $G'$ may be regarded as a plane graph with one ring $C$,
and as such it is $C$-critical.
\end{lemma}

\proof If $G'$ is not $C$-critical, then let $e\in E(G')\setminus E(C)$ be an edge such that every precoloring of $C$
that extends to $G'-e$ also extends to $G'$.  Observe that every precoloring of ${\cal R}$ that extends to $G-e$ also extends
to $G$, contrary to the assumption that $G$ is $\cal R$-critical.~\qed
\bigskip

Critical plane graphs with one ring of length at most twelve were described by Thomassen~\cite{thom-surf} and independently
by Walls~\cite{walls-enum} (actually, both papers describe $\phi$-critical graphs for some fixed precoloring $\phi$
of $R$, but Theorem~\ref{thm-planechar} follows straightforwardly from the characterizations):

\begin{theorem}\label{thm-planechar}
Let $G$ be a plane graph of girth $5$ with one ring $R$ such that 
$|V(R)|\le 12$.  If $G$ is $R$-critical and
$R$ is an induced cycle, then

\begin{itemize}
\item[(a)] $|V(R)|\ge 9$ and $G-V(R)$ is a tree with at most $|V(R)|-8$ vertices, or
\item[(b)] $|V(R)|\ge 10$ and $G-V(R)$ is a connected graph with at most 
$|V(R)|-5$ vertices
that contains exactly one cycle, and the length of this cycle is $5$, or
\item[(c)] $|V(R)|=12$ and every second vertex of $R$ has degree two and is contained in a facial $5$-cycle.
\end{itemize}
\end{theorem}

In this section, we generalize this result by giving bounds on the weight of planar critical graphs
with one ring.

\begin{theorem}
\mylabel{thm:diskgirth5}
Let $\epsilon\le 1/1278$ be a positive
real number and let $s:\{5,6,\ldots\}\to{\mathbb R}$ be an increasing function
satisfying conditions (S1)--(S5) formulated in Sections~\ref{sec-disch} and \ref{sec-winners}.
Let $G$ be a plane graph with one ring $R$ of length $l\ge 5$
such that $G$ is $R$-critical
and has no cycle of length at most four, and let $w$ be the weight function arising from $s$ as described in Section~\ref{sec-disch}.  Then
\begin{itemize}
\item $l\ge 8$ and $w(G,\{R\})\le s(l-3)+s(5)$, and furthermore,
\item if $R$ does not satisfy (E1), then $l\ge 9$ and $w(G,\{R\})\le s(l-4)+2s(5)$,
\item if $(G,R)$ is not very exceptional, then $l\ge 10$ and $w(G,\{R\})\le s(l-5)+5s(5)$, and
\item if $(G,R)$ is not exceptional, then $l\ge 11$ and $w(G,\{R\})\le s(l-5)-5s(5)$.
\end{itemize}
\end{theorem}

\proof
Let us note that $s(l-4)+2s(5)\le s(l-3)+s(5)$ for $l\ge9$ by (S4),
and hence whenever $G$ satisfies the second conclusion, it satisfies the first.
If $G$ satisfies $(E1)$, then $l\ge 8$ and $G$ has a face of length $a$ such that
$a\le 7$.  We can assume that the other face of $G$ is at least as long as $a$,
that is, $l+2-a\ge a$.  Then, $w(G,\{R\})=s(a)+s(l+2-a)\le s(l-3)+s(5)$,
where the inequality holds by convexity.
If $G$ satisfies (E2), then it is very exceptional, $l\ge9$ and
$w(G,\{R\})=s(l-4)+2s(5)$.  
If $G$ satisfies (E3), then it is very exceptional,
$l\ge 11$ and $w(G,\{R\})=s(l-5)+s(5)+s(6)\le s(l-4)+2s(5)$, where the inequality
follows from convexity.
If $G$ satisfies (E4) or (E5), then $l\ge 10$
and $w(G,\{R\})\le s(l-5)+5s(5)\le s(l-4)+2s(5)$, where the second inequality
follows from convexity and (S4).  Finally, suppose that $G$ is not exceptional.
By Theorem~\ref{thm-planechar}, we have $l\ge 11$; thus
$s(l-5)-5s(5)\le s(l-5)+5s(5)\le s(l-4)+2s(5)$ by convexity and (S4).
Therefore, it suffices to prove that $w(G,\{R\})\le s(l-5)-5s(5)$
whenever $G$ is not exceptional.

Suppose for a contradiction that $(G,R)$ is not exceptional, 
and yet $w(G,\{R\})> s(l-5)-5s(5)$.
We may assume that the theorem holds for all graphs with fewer edges than $G$.

\claimno=0

{
\newclaim{lowbd}{$l\ge12$}
}

\noindent
To prove (\refclaim{lowbd}) let $l\le11$. Since $G$ is not exceptional,
it follows from Theorem~\ref{thm-planechar} that $l=11$, every
face of $G$ has length five, and there are at most seven
faces. Thus $w(G,\{R\})\le 7w(5)\le s(6)-5w(5)=s(l-5)-5w(5)$ by (S4), 
a contradiction.
This proves (\refclaim{lowbd}).

{
\newclaim{cl1}{There is no path of length at most two with both ends in $R$ that
is otherwise disjoint from $R$ (i.e., $G$ satisfies {\rm (I4)}).}
}

\noindent
To prove (\refclaim{cl1}) let $P$ be a path in $G$ of length one or two with ends 
$u,v\in V(R)$, and otherwise disjoint from $R$.
Let $C_1,C_2$ be the two cycles of $R\cup P$ other than $R$,
and for $i=1,2$ let $G_i$ be the subgraph of $G$ drawn in the closed disk bounded by $C_i$ and $l_i=|C_i|$.
Note that $l_1+l_2=l+2|\rt{E(}P)|$.

Since $G$ does not satisfy (E1) and satisfies (I0), we can assume
that $G_1\neq C_1$.  Hence $G_1$ is $C_1$-critical by Lemma~\ref{lem:diskcritical}.
Assume for a moment that $G_2=C_2$. 
If $G_1$ is not very exceptional, then using the minimality of $G$, we have
$w(G,\{R\})=w(G_1,\{C_1\})+s(l_2)\le s(l_1-5)+5s(5)+s(l_2)\le s(l_1+l_2-10)+6s(5)\le s(l-5)-5s(5)$
by the convexity and (S4), a contradiction.  Similarly, we exclude the case that $P$ has length one
and $G_1$ is very exceptional.  
Finally, if $G_1$ is very exceptional and $|\rt{E(}P)|=2$,
then $G\rt{-} V(R)$ consists of one or two 
adjacent vertices of degree three in $G$.  Let $a_1\le a_2\le \ldots$ be the lengths of the faces of $G$.
If \rt{$G-V(R)$ has one vertex}, then since $G$ does not satisfy 
(E2) \rt{or} (E3), we have \rt{either} $a_1\ge 6$ (and $l\ge 12$),
in which case $w(G,\{R\})=s(a_1)+s(a_2)+s(a_3)\le 2s(6)+s(l-6)\le s(l-5)-5s(5)$, by convexity and (S4),
or $a_1=5$ and $a_2\ge 7$ (and $l\ge 13$), 
in which case
$w(G,\{R\})=s(a_1)+s(a_2)+s(a_3)\le s(5)+s(7)+s(l-6)\le s(l-5)-5s(5)$, 
again by convexity and (S4).
If \rt{$G-V(R)$ consists of two adjacent vertices of degree three,} then,
 since $G$ does not satisfy (E4), we have $a_3\ge 6$ and $l\ge 12$; thus
$w(G,\{R\})=s(a_1)+s(a_2)+s(a_3)+s(a_4)\le 2s(5)+s(6)+s(l-6)\le s(l-5)-5s(5)$.
This is a contradiction.

Thus we may assume that $G_1\ne C_1$ and $G_2\ne C_2$.
Therefore, $G_1$ is $C_1$-critical and $G_2$ is $C_2$-critical by 
Lemma~\ref{lem:diskcritical}.
Furthermore, we may assume that $P$ cannot be chosen so that $G_2=C_2$.
That implies that $G_1$ and $G_2$ are not very exceptional, and hence
$w(G,\{R\})\le s(l_1-5)+5s(5)+s(l_2-5)+5s(5)\le s(l-5)-5s(5)$.
a contradiction. This proves (\refclaim{cl1}).
\bigskip

Let $\phi$ be a precoloring of $R$ that does not extend to a $3$-coloring of $G$.

{
\newclaim{cl0}{$G$ is $\phi$-critical.}
}

\noindent To prove (\refclaim{cl0}) suppose to the contrary
that $G$ is not $\phi$-critical.
Then $G$ contains a proper $\phi$-critical subgraph $G'$.  
By Lemma~\ref{lem:crit3conn}, $G'$ is $2$-connected, and thus all its faces are bounded by cycles.  
Note that $G'$ is not very exceptional
by (\refclaim{cl1}). Since $G'$ has fewer edges than $G$, we have $w(G',\{R\})\le s(l-5)+5s(5)$ by induction.
For $f\in{\cal F}(G')$ let $G_f$ be the subgraph of $G$ drawn inside the 
closure of $f$, and let $C$ be the cycle bounding $f$.
By Lemma~\ref{lem:diskcritical}, $G_f$ is either equal to $C$,
or it is $C$-critical. 
Thus by induction, the convexity of $s$ and (S4), we have 
$w(G_f,\{C\})\le s(|f|)$.
Furthermore, if $G_f$ is not equal to $C$, then $w(G_f,\{C\})\le s(|f|-3)+s(5)$.
Let $f_0$ be a face of $G'$ such that $G_{f_0}$ is not equal to $f_0$.  Note that $|f_0|\ge 8$
by Theorem~\ref{thm-planechar}.  We have
\begin{eqnarray*}
w(G,\{R\})&=&\sum_{f\in{\cal F}(G)} s(|f|) =
\sum_{f'\in{\cal F}(G')} w(G_{f'},\{f'\})\\
&\le&s(|f_0|-3)+s(5)-s(|f_0|)+\sum_{f'\in{\cal F}(G')}s(|f'|)\\
&=&s(|f_0|-3)+s(5)-s(|f_0|)+w(G',\{R\})\\
&\le&s(|f_0|-3)-s(|f_0|)+s(l-5)+6s(5)\le s(l-5)-5s(5),
\end{eqnarray*}
where the last inequality holds by convexity and (S4).
This proves (\refclaim{cl0}).

{
\newclaim{cl2}{The graph $G$ does not have two adjacent vertices of degree two
(i.e., $G$ satisfies {\rm (I5)}). Furthermore, every vertex of degree two is incident
with a face of length at most six.}
}

\noindent
To prove (\refclaim{cl2}) let $u$ and $v$ be two adjacent vertices of degree two in $R$.
Let $G'$ and $R'$ be the graphs obtained from $G$ and $R$, respectively,
by identifying $u$ and $v$ into a single vertex $w$.  Let $\phi'$ be a $3$-coloring of $R'$ matching $\phi$
on $R'-w$.  Note that $G'$ is $\phi'$-critical, and $G'$ contains no cycle of length at most four by (\refclaim{cl1}).
Let $d$ be the length of the face $f$ of $G$ incident with the edge $uv$.
By (\refclaim{cl1}), if $G'$ is exceptional, then it satisfies (E5); hence
$G$ has four faces of length five, a $6$-face and a face of length $l-6$
and $w(G,\{R\})=s(l-6)+s(6)+4s(5)\le s(l-5)-5s(5)$
by (\refclaim{lowbd}) and (S4).  
Therefore, assume that $G'$ is not exceptional.  By the minimality of $G$ we have $w(G',\{R'\})\le s(l-6)-5s(5)$,
and since the face corresponding to $f$ contributes $s(d-1)$ to $w(G',\{R'\})$, we conclude that $d-1<l-6$. Thus
$w(G,\{R\})=w(G',\{R'\})-s(d-1)+s(d)\le s(l-6)-5s(5)-s(d-1)+s(d)\le
s(l-5)-5s(5)$ by convexity.  The case that a vertex $v$ of degree two is incident with a face
of length at least $7$ is handled similarly, with $G'$ obtained either by suppressing $v$ or by
identifying its neighbors, depending on whether the colors of these neighbors according to $\phi$
differ or not.  This proves (\refclaim{cl2}).
\medskip

{
\newclaim{cl3}{A good configuration appears in $G$.}
}

\noindent
To prove (\refclaim{cl3}) suppose for a contradiction that no good configuration
appears in $G$.
By Lemma~\ref{lem:i012} the graph satisfies (I0), (I1) and (I2).
By Lemma~\ref{lem:crit3conn}, the graph $G$ satisfies (I3) and (I6).
By (\refclaim{cl1}) and (\refclaim{cl2}) it satisfies (I4) and (I5).
The assumptions (I7) and (I8) are trivially satisfied by planar graphs with only one ring.
Let $M$ be the null graph.
We deduce from Lemma~\ref{lem:noconfigweight} that
 $w(G,{\cal R})\le 4n_2/3+52\epsilon n_3-8$,
where $n_2$ and $n_3$ are as in Lemma~\ref{lem:fincharge}.
By (I5) we have $n_2\le l/2$, thus $4n_2/3+52\epsilon n_3\le (2/3+26\epsilon)l$.
If $l\ge16$, then 
$$w(G,{\cal R})\le (2/3+26\epsilon)l-8\le l-13-10\epsilon=s(l-5)-5s(5)$$
because $\epsilon\le1/1278$, a contradiction.
Thus we may assume that $l\le 15$, and hence $n_2\le7$.
If $l=15$, then $w(G,{\cal R})\le 28/3+8\cdot 52\epsilon-8\le l-13-10\epsilon=s(l-5)-5s(5)$,
again a contradiction.
If $l=13$, then we $n_2\le 6$ and $w(G,{\cal R})\le 7\cdot 52\epsilon\le s(8)-5s(5)$.

Suppose that $l=12$.
If $n_2\le 5$, then 
$w(G,{\cal R})\le 20/3+52\cdot12\epsilon-8\le0\le s(7)-5s(5)$,
because $270\epsilon\le s(7)\le s(8)/3\le(s(7)+s(9))/6$ by (S4), implying
that $\epsilon\le1/1350$. Thus we may assume that $n_2=6$ and $n_3=6$.
By Theorem~\ref{thm-planechar}, all faces sharing an edge with $R$
have length $5$, thus the internal vertices that have a neighbor in $R$ form a $6$-cycle $K$.
By Lemma~\ref{lem:diskcritical} and Theorem~\ref{thm-planechar}, we have that $K$ bounds a face,
thus all its vertices have degree three.  This contradicts (I1).  It follows that
if $l=12$ and $n_2=6$, then $n_3\le 5$; thus $w(G,{\cal R})\le 260\epsilon \le s(7)-5s(5)$ by (S4).

Thus by (\refclaim{lowbd}) we may assume that $l=14$.  
If $n_2\le6$, then we have $w(G,{\cal R})\le 8\cdot 52\epsilon\le s(9)-5s(5)$;
hence $n_2=7$.  Furthermore, using Lemma~\ref{lem:newconfigweight}, we conclude that $b=0$, where
$b$ is as in that lemma. 
Then vertices of degree two and three alternate on $R$, and every 
face that shares an edge with $R$ has length five. The neighbors of the vertices of $R$
of degree three form a $7$-cycle, which bounds a face by Theorem~\ref{thm-planechar}.
Then, $w(G,\{R\})=s(7)+7s(5)\le s(9)-5s(5)$.  This proves (\refclaim{cl3}).
\medskip

{
\newclaim{cl4}{The graph $G$ is well-behaved.}
}

\noindent
To prove (\refclaim{cl4}), assume to the contrary that $G$ is not well-behaved.
Thus there exists a path $P$ of length at most four, with ends $u,v\in V(R)$ and otherwise
disjoint from $R$, that is not allowable.
We may assume that $P$ is such a path of the shortest
possible length. By (\refclaim{cl1}), the path $P$ has length at least three.

Let $C_1$, $C_2$, and $R$ be the three cycles of $R\cup P$, and for $i=1,2$
let $G_i$ be the subgraph of $G$ consisting of all vertices and edges
drawn in the closed disk bounded by $C_i$.
We claim that $C_1$ and $C_2$ are induced cycles.
To prove this claim suppose to the contrary that some edge has ends
$x,y\in V(C_i)$ for some $i\in\{1,2\}$, but that the edge itself does not
belong to $C_i$.
Then one of the vertices $x$ and $y$, say $x$, belongs to the interior of $P$, and
$y$ does not belong to $P$.
By (\refclaim{cl1}), the vertex $x$ is not a neighbor of $u$ or $v$, and hence
$P$ has length four, and $x$ is the middle vertex of $P$.
Let the vertices of $P$ be $u,u',x,v',v$, in order.
Since $P$ was chosen minimal, the two paths $uu'xy$ and $vv'xy$
are allowable; hence $G_i$ consists of two $5$-faces and the path $P$ is allowable,
a contradiction.  This proves that $C_1$ and $C_2$ are induced cycles.

It follows from (\refclaim{cl1}) and (\refclaim{cl2}) that $G_1$ and $G_2$ are not very exceptional
and that $G_i\ne C_i$. By Lemma~\ref{lem:diskcritical} the graph
$G_i$ is $C_i$-critical for $i=1,2$.
Let $l_i=|C_i|$.
By induction we have
\begin{eqnarray*}
w(G,\{R\}) &=& w(G_1,\{C_1\})+ w(G_2,\{C_2\})\\
&\le& s(l_1-5)+5s(5)+s(l_2-5)+5s(5)\\
&\le& s(l_1+l_2-15)+11s(5)\le s(l-5)-5s(5),
\end{eqnarray*}
by convexity and (S4). This proves (\refclaim{cl4}).
\bigskip

It follows from (\refclaim{cl3}), (\refclaim{cl4}) and Lemma~\ref{lem:strongappear} 
that some good configuration strongly appears in $G$, for if the
second outcome of Lemma~\ref{lem:strongappear} holds, then
$(G,R)$ either is exceptional or satisfies the conclusion of the theorem.
Let $\gamma$ be a good configuration that strongly appears in $G$,
and let $G'$ be the $\gamma$-reduction of $G$.
By Lemma~\ref{lemma-color} the $3$-coloring $\phi$ does not extend to
a $3$-coloring of $G'$.
Thus $G'$ has a $\phi$-critical subgraph $G''$.
By Lemma~\ref{lem:4cycles} the graph $G''$ has no cycles of length at most four
($G$ satisfies (I9) by Lemma~\ref{lem:diskcritical} and Theorem~\ref{thm-planechar}).
By Lemma~\ref{lem:crit3conn}, the graph $G''$ satisfies (I3) and (I6).
Since $G$ is $\phi$-critical by (\refclaim{cl0}), 
$G''$ is not a subgraph of $G$; hence $G''$ contains a new vertex or edge.

For a face $f''$ of $G''$ let  
$G_{f''}^1$, $G_{f''}^2$, \ldots, $G_{f''}^{k_{f''}}$
be the members of the $G$-expansion of $S_{f''}$,
defined in Section~\ref{sec-winners},
and let $C_{f''}^1$, $C_{f''}^2$, \ldots, $C_{f''}^{k_{f''}}$
be the corresponding rings so that $C_{f''}^i$ is a subgraph of $G_{f''}^i$.

{
\newclaim{cl5}{Let $f''$ be a face of $G''$. Then
$$\sum_{i=1}^{k_{f''}} w(G_{f''}^i,\{C_{f''}^i\})\le s(|f''|)-c(f'').$$}
}

Note that by Lemma~\ref{lem:diskcritical}, we have that either $G_{f''}^i=C_{f''}^i$ or $G_{f''}^i$ is $C_{f''}^i$-critical
for each $i$.  To prove (\refclaim{cl5}), let us discuss the possible cases in the definition of the contribution of a face:
\begin{itemize}
\item If  $|S_{f''}|=1$ and
$G^1_{f''}$ satisfies (E0), then by Lemma~\ref{lem:winners} we have $c(f'')\neq-\infty$, hence $f''$ has zero elasticity,
$c(f'')=0$ and $w(G_{f''}^1,\{C_{f''}^1\})=s(|f''|)$.
\item If  $|S_{f''}|=1$ and
$G^1_{f''}$ satisfies (E1), then similarly we have $\el(f'')<5$ and $c(f'')=s(8-\el(f''))-2s(5)$.  Note that by
Lemma~\ref{lem:elasticity} we have $\el(f'')\le 3$.  By induction, $w(G_{f''}^1,\{C_{f''}^1\})\le s(|C_{f''}^1|-3)+s(5)=s(|f''|+\el(f'')-3)+s(5)$,
and $s(|f''|+\el(f'')-3)+s(5)\le s(|f''|)-s(8-\el(f''))+2s(5)$ by convexity.
\item If  $|S_{f''}|=1$ and
$G^1_{f''}$ satisfies (E2), then $\el(f)\le 3$, $c(f'')=s(9-\el(f''))-3s(5)$,
and $w(G_{f''}^1,\{C_{f''}^1\})=s(|C_{f''}^1|-4)+2s(5)=s(|f''|+\el(f'')-4)+2s(5)\le s(|f''|)-c(f'')$ by convexity.
\item If  $|S_{f''}|=1$ and
$G^1_{f''}$ satisfies (E3), then $c(f'')=s(11-\el(f''))-2s(6)-s(5)$ and
$w(G_{f''}^1,\{C_{f''}^1\})=s(|C_{f''}^1|-5)+s(5)+s(6)=s(|f''|+\el(f'')-5)+s(5)+s(6)\le s(|f''|)-c(f'')$.
\item If  $|S_{f''}|=1$ and
$G^1_{f''}$ satisfies (E4) or (E5), then $c(f'')=s(10-\el(f''))-6s(5)$ and
$w(G_{f''}^1,\{C_{f''}^1\})\le s(|C_{f''}^1|-5)+5s(5)=s(|f''|+\el(f'')-5)+5s(5)\le s(|f''|)-c(f'')$.
\item Suppose that $k_{f''}=2$, $G_{f''}^1=C_{f''}^1$ and $G_{f''}^2=C_{f''}^2$, where $|C_{f''}^1|\le |C_{f''}^2|$.  If $|C_{f''}^1|=5$,
then $c(f'')=s(10-\el(f''))-6s(5)$ and $w(G_{f''}^1,\{C_{f''}^1\})+w(G_{f''}^2,\{C_{f''}^2\})=s(|C_{f''}^2|)+s(5)=s(|f''|+\el(f'')-5)+s(5)<s(|f''|)-c(f'')$.
Otherwise, $c(f'')=s(12-\el(f''))-2s(6)$ and
$w(G_{f''}^1,\{C_{f''}^1\})+w(G_{f''}^2,\{C_{f''}^2\})=s(|C_{f''}^1|)+s(|C_{f''}^2|)\le s(6)+s(|f''|+\el(f'')-6)\le s(|f''|)-c(f'')$.
\item Suppose that $k_{f''}=1$ and $G^1_{f''}$ is not exceptional.
\begin{itemize}
\item Let us consider the case that $G_{f''}^1$ contains a path $P=v_1v_2v_3v_4$ such that $v_1,v_4\in V(C_{f''}^1)$, $v_2,v_3\not\in V(C_{f''}^1)$ and both of the open disks
$\Delta_1$ and $\Delta_2$ bounded by $P$ and paths in $C_{f''}^1$ contain at least two vertices of $G$.  In this case, $c(f'')=s(7)$.  Let $H_i$ be the
subgraph of $G_{f''}^1$ drawn in $\Delta_i$ and $K_i$ the cycle bounding $\Delta_i$, for $i\in\{1,2\}$.  Neither $H_1$ nor $H_2$ is very exceptional,
thus we have $w(G_{f''}^1,\{C_{f''}^1\})=w(H_1,K_1)+w(H_2,K_2)\le s(|K_1|-5)+s(|K_2|-5)+10s(5)\le
s(|K_1|+|K_2|-15)+11s(5)=s(|f''|+\el(f'')-9)+11s(5)<s(|f''|)-s(7)$, 
since $\el(f'')\le 5$ and $|f''|+\el(f'')\ge |K_1|+|K_2|-6\ge 14$,
because $G''$ satisfies (I9) by Lemma~\ref{lem:diskcritical} 
and Theorem~\ref{thm-planechar}.
\item Otherwise, $c(f'')=s(11-\el(f''))-s(6)+5s(5)$.  In this case, we have
$w(G_{f''}^1,\{C_{f''}^1\})\le s(|C_{f''}^1|-5)-5s(5)=s(|f''|+\el(f'')-5)-5s(5)\le s(|f''|)-c(f'')$.
\end{itemize}
\item If $k_{f''}=2$ and $G_{f''}^1\neq C_{f''}^1$, then $c(f'')=s(12-\el(f''))-2s(6)$ and
$w(G_{f''}^1,\{C_{f''}^1\})+w(G_{f''}^2,\{C_{f''}^2\})\le
s(|C_{f''}^1|-3)+s(5)+s(|C_{f''}^2|)\le s(|f''|+\el(f'')-8)+2s(5)<s(|f''|)-c(f'')$
\item If $k_{f''}\ge 3$, then $c(f'')=s(12-\el(f''))-2s(6)$
and $\sum_{i=1}^{k_{f''}} w(G_{f''}^i,\{C_{f''}^i\})\le s(|f''|+\el(f'')-(k_{f''}-1)5)+(k_{f''}-1)s(5)<s(|f''|)-c(f'')$.
\end{itemize}
Therefore, in all the cases, (\refclaim{cl5}) holds.
\bigskip

By Lemma~\ref{lem:facecover}, we have $w(G,\{R\})\le \delta+\sum_{f''\in {\cal F}(G'')}\sum_{i=1}^{k_{f''}}w(G_{f''}^i,\{C_{f''}^i\})$,
where $\delta=s(6)$ if $\gamma$ is isomorphic to $\R3$ and $\delta=0$ otherwise.
By (\refclaim{cl5}) this implies that
\begin{eqnarray*}
w(G,\{R\})&\le& \delta+\sum_{f''\in {\cal F}(G'')} (s(|f''|)-c(f''))\\
&=&w(G'',\{R\})+\delta-\sum_{f''\in {\cal F}(G'')}c(f'')\\
&=&w(G'',\{R\})-c(G'').
\end{eqnarray*}
By Lemma~\ref{lemma-rednt}, $G''$ is not very exceptional; 
hence $w(G'',\{R\})\le s(l-5)+5s(5)$ by induction.
Note that $c(G'')\ge 10s(5)$ by Lemma~\ref{lem:winners}; thus
$$w(G,\{R\})\le w(G'',\{R\})-c(G'')\le s(l-5)-5s(5),$$
which is a contradiction finishing the proof of the theorem.~\qed

\bigskip

\noindent
{\bf Proof of Theorem~\ref{thm:main}.}
Let 
\rt{%
$\epsilon=2/4113$, 
$s(5)= 4/4113$, 
$s(6)=72/4113$, 
$s(7)=540/4113$,
$s(8)=2184/4113$,
}
and $s(l)=l-8$ for $l\ge9$.
Then conditions (S1)--(S4) hold. Furthermore, note that $s(l-3)+s(5)<s(l)$ for all $l\ge 8$.

Let $G$ be a plane graph of girth at least five
with a cycle $C$ and let $\phi$ be a precoloring of $C$ that 
does not extend to a $3$-coloring of $G$.  
We may assume that $G$ is $\phi$-critical, and hence $C$ is a face of $G$.  
By Theorem~\ref{thm:diskgirth5}, we have $w(G,\{C\})\le s(|C|-3)+s(5)<|V(C)|$.  
Note that $3|V(G)|-2|V(C)|=\sum_f |f|\le \sum_f 5s(|f|)/s(5)=5w(G,\{C\})/s(5)$, 
where the sum is over all faces of $G$, except the one bounded by $C$.
Therefore, $|V(G)|\le (5/s(5)+2)|V(C)|/3\le 1715|V(C)|$, as desired.~\qed

\section{Summary}\label{sec-summary}

In this section, we provide a summary result that will be used as a basis for the proofs in the next
paper of the series, to avoid the need to repeat many of the definitions used here.
Let $\Pi$ be a surface with boundary and $c$ a simple curve intersecting the boundary of $\Pi$ exactly in its ends.
The compact topological space obtained from $\Pi$ by cutting along $c$ (i.e., removing $c$ and adding two new pieces of boundary
corresponding to $c$) is a union of at most two surfaces.
If surfaces $\Pi_1,\ldots, \Pi_k$ are obtained from $\Pi$ by repeating this construction,
we say that they are \emph{fragments} of $\Pi$.

Consider a graph $H$ embedded in $\Pi$ with rings $\QQ$, and let $f$ be a face of $H$.
Let $\Pi_f$ be a surface whose interior is homeomorphic; note that $\Pi_f$ is unique up to homeomorphism
and that the cuffs of $\Pi_f$ correspond to the facial walks of $f$.

Let $G$ and $G'$ be $\RR$-critical graphs embedded in $\Sigma$ with rings $\RR$.
Suppose that there exists a collection $\{(J_f,S_f):f\in F(G')\}$ of subgraphs $J_f$ of $G$ and sets $S_f$ of faces of $J_f\cup\bigcup\RR$ 
such that $J_f$ is the union of the boundary walks of $S_f$, and a set $X\subset F(G)$ such that 
\begin{itemize}
\item for every $f\in F(G')$, the boundary of $S_f$ is not equal to the union of $\RR$,
\item for every $f\in F(G')$, the surfaces of the $G$-expansion of $S_f$ are fragments of $\Sigma_f$,
\item for every face $h\in F(G)\setminus X$, there exists unique $f\in F(G')$ such that $h$ is subset of a member of $S_f$, and
\item if $X\neq \emptyset$, then $X$ consists of a single closed $2$-cell face of length $6$.
\end{itemize}
We say that $X$ together with this collection forms a \emph{cover of $G$ by faces of $G'$}.
For a face $f\in F(G')$, let $\el(f)=\left(\sum_{h\in S_f} |h|\right)-|f|$.

\begin{theorem}
Let $G$ be a well-behaved graph embedded in a surface $\Sigma$ with rings $\RR$
satisfying {\rm (I0)--(I9)} and let $M$ be a subgraph of $G$ that captures $(\le\!4)$-cycles.
Let $\ell(\RR)$ be the sum of the lengths of the rings in $\RR$ and $g$ the Euler genus of $\Sigma$,
and assume that $g>0$ or $|\RR|>1$.
Let $\epsilon$ be a real number satisfying $0<\epsilon\le 1/1278$,
let $s:\{5,6,\ldots\}\to{\mathbb R}$ be a function satisfying {\rm (S1)--(S4)},
and suppose that $w(G,\RR)>8g+8|{\cal R}|+(2/3+26\epsilon)\ell(\RR)+20|E(M)|/3-16$.
If $G$ is $\RR$-critical, then there exists an $\RR$-critical graph $G'$ embedded in $\Sigma$ with rings $\RR$
such that $|E(G')|<|E(G)|$, every vertex-like ring of $G$ is also vertex-like in $G'$, and the following conditions hold.
\begin{enumerate}
\item If $G$ has girth at least five, then there exists a set $Y\subseteq V(G')$ of size at most two such that $G'-Y$ has
girth at least five.  
\item If $C'$ is a $(\le 4)$-cycle in $G'$, then $C'$ is non-contractible and $G$ contains a non-contractible
cycle $C$ of length at most $|C'|+3$ such that $C\not\subseteq M$.  Furthermore, all ring vertices of $C'$ belong to $C$,
and if $C'$ is a triangle disjoint from the rings
and its vertices have distinct pairwise non-adjacent neighbors in a ring $R$ of length $6$, then
$G$ contains edges $cr$ and $c'r'$ with $c,c'\in V(C)\setminus V(R)$ and $r,r'\in V(R)\setminus V(C)$
such that $r$ and $r'$ are non-adjacent.
\item $G'$ has a face that either is not semi-closed $2$-cell or has length at least $6$.
\item There exists $X\subset F(G)$ and a collection $\{(J_f,S_f):f\in F(G')\}$ forming a cover of $G$ by faces of $G'$ satisfying
the following conditions.
\begin{enumerate}
\item For a semi-closed $2$-cell or omnipresent face $f\in F(G')$, let its contribution
$c(f)$ be defined as in Section~\ref{sec-winners}.  Then $\sum_{f\in F(G')} \el(f)\le 10$ and if $f$ is an omnipresent face, then $\el(f)\le 5$.
Furthermore, if every face of $G'$ is semi-closed $2$-cell or omnipresent, $G'$ satisfies (I6), and every vertex-like ring of $G'$ is also vertex-like in $G$, then $\sum_{f\in F(G')} c(f)\ge |X|s(6)$. 
\item If every vertex-like ring of $G'$ is also vertex-like in $G$, $f\in F(G')$ is semi-closed $2$-cell and $G_1$, \ldots, $G_k$ are the components of the $G$-expansion of $S_f$, where for $1\le i\le k$, $G_i$ is embedded in
a disk with one ring $R_i$, then $\sum_{i=1}^k w(G_i,\{R_i\})\le s(|f|)-c(f)$.
\end{enumerate}
\end{enumerate}
\end{theorem}
\begin{proof}
Let $n_2$ be the number of ring vertices of degree two not belonging to $M$ and let $n_3$ be the number of ring vertices of degree three.
By (I5) we have $n_2\le \ell(R)/2$, and since $n_2+n_3\le \ell(R)$, we have
$4n_2/3+52\epsilon n_3\le (2/3+26\epsilon)\ell(R)$.  Consequently, 
$w(G,\RR)>8g+8|{\cal R}|+4n_2/3+52\epsilon n_3+20|E(M)|/3-16$, and
by Lemma~\ref{lem:noconfigweight}, a good configuration $\gamma$ appears in $G$ and does not touch $M$.  By Lemma~\ref{lem:strongappear},
we can assume that $\gamma$ appears strongly in $G$.
Let $\phi$ be a precoloring of $\RR$ that does not extend to a $3$-coloring of $G$, and let $G_1$ be a $\gamma$-reduction
of $G$ with respect to $\phi$.  By Lemma~\ref{lemma-color}, $\phi$ does not extend to a $3$-coloring of $G_1$,
and thus $G_1$ contains an $\RR$-critical subgraph $G'$.
Since $G'$ is $\RR$-critical, it satisfies (I0). Clearly, $|E(G')|<|E(G)|$.  Let us now show that $G'$ has the required properties:
\begin{enumerate}
\item If $G$ has girth at least five, then every $(\le\!4)$-cycle in $G'$ contains a new vertex or a new edge, and thus they can all be intersected by at most two vertices.
\item Follows from Lemma~\ref{lem:4cycles}.
\item Suppose that all faces of $G'$ are semi-closed $2$-cell.  In particular, $G'$ does not have an omnipresent face,
and thus it satisfies (I6).  If $G'$ contains a new edge or a new vertex, then the claim holds by Lemma~\ref{lem:winners}.
Otherwise, $G'$ is a proper subgraph of $G$, and thus there exists a cycle $C$ bounding a face in $G'$, but not in $G$.
Let $H$ be the subgraph of $G$ drawn in the closed disk bounded by $C$.  By Lemma~\ref{lem:diskcritical},
$H$ is $C$-critical, and by Theorem~\ref{thm-planechar}, we conclude that $C$ has length at least $8$.
\item For each $f\in F(G')$, we define $S_f$ and $J_f$ as in Section~\ref{sec-winners}.  As $G'$ is not equal to the union of $\RR$,
the boundary of $S_f$ is distinct from the union of $\RR$ for each $f\in F(G')$.  
The construction of $J_f$ and $S_f$ ensures that the surfaces corresponding to the faces of $S_f$ are constructed
from $\Sigma_f$ by cutting along simple curves with ends in cuffs, as described in the definition of fragments.
By Lemma~\ref{lem:facecover}, there exists a set $X\subset F(G)$ (which is either empty or consists of one $6$-face)
such that $X$ together with the collection $\{(J_f,S_f):f\in F(G')\}$ is a cover of $G$ by faces of $G'$.  Let us argue that it
has the properties asserted by the lemma:
\begin{enumerate}
\item The first part follows from Lemma~\ref{lem:elasticity}.  If $G'$ contains a new vertex or a new edge, then the second
part follows from Lemma~\ref{lem:winners}.  Otherwise, $G'$ is a proper subgraph of $G$ and all its faces have elasticity $0$.
If $f$ is a semi-closed $2$-cell of $G'$, then $c(f)\ge 0$, and if additionally $f$ is not a face of $G$, then $c(f)\ge s(8)-2s(5)>s(6)$.
If $f$ is an omnipresent face, then note that no component of $G'$ satisfies (E1), (E2) or (E3), since $G$ satisfies (I4).
Since $G'$ is $\RR$-critical, at least one component of $G'$ does not satisfy (E0), and thus $c(f)\ge 5-5s(5)>s(6)$.
Since $G'\neq G$, we conclude that $\sum_{f\in F(G')} c(f)>s(6)\ge |X|s(6)$.
\item This was proved as (\refclaim{cl5}) in Section~\ref{sec-disk}.
\end{enumerate}
\end{enumerate}
\end{proof}

\bibliographystyle{acm}
\bibliography{4critsurf}

\begin{thebibliography}{10}

\bibitem{BonMur}
{\sc Bondy, J., and Murty, U.}
\newblock {\em Graph {T}heory with {A}pplications}.
\newblock North-Holland, New York, Amsterdam, Oxford, 1976.

\bibitem{dvkaw}
{\sc Dvo\v{r}\'ak, Z., and Kawarabayashi, K.}
\newblock {Choosability of planar graphs of girth 5}.
\newblock {\em ArXiv e-prints 1109.2976\/} (2011).

\bibitem{DvoKawTho}
{\sc Dvo\v{r}\'ak, Z., Kawarabayashi, K., and Thomas, R.}
\newblock Three-coloring triangle-free planar graphs in linear time.
\newblock {\em Trans. on Algorithms 7\/} (2011), article no. 41.

\bibitem{proof-part3}
{\sc Dvo\v{r}\'ak, Z., Kr\'al', D., and Thomas, R.}
\newblock Three-coloring triangle-free graphs on surfaces {III}. {G}raphs of
  girth five.
\newblock {\em ArXiv e-prints 1402.4710\/} (Feb. 2014).

\bibitem{proof-part4}
{\sc Dvo\v{r}\'ak, Z., Kr\'al', D., and Thomas, R.}
\newblock Three-coloring triangle-free graphs on surfaces {IV}. {B}ounding face
  sizes of $4$-critical graphs.
\newblock {\em ArXiv e-prints 1404.6356v4\/} (May 2015).

\bibitem{proof-havel}
{\sc Dvo\v{r}\'ak, Z., Kr\'al', D., and Thomas, R.}
\newblock Three-coloring triangle-free graphs on surfaces {V}. {C}oloring
  planar graphs with distant anomalies.
\newblock {\em ArXiv e-prints 0911.0885v3\/} (Jan. 2016).

\bibitem{gimbel}
{\sc Gimbel, J., and Thomassen, C.}
\newblock Coloring graphs with fixed genus and girth.
\newblock {\em Trans. Amer. Math. Soc. 349\/} (1997), 4555--4564.

\bibitem{Gro}
{\sc Gr\"otzsch, H.}
\newblock Ein {D}reifarbensatz f\"ur dreikreisfreie {N}etze auf der {K}ugel.
\newblock {\em Wiss. Z. Martin-Luther-Univ. Halle-Wittenberg Math.-Natur. Reihe
  8\/} (1959), 109--120.

\bibitem{conj-havel}
{\sc Havel, I.}
\newblock On a conjecture of {G}r{\"u}nbaum.
\newblock {\em J. Combin. Theory Ser. B 7\/} (1969), 184--186.

\bibitem{koyan}
{\sc Kostochka, A.~V., and Yancey, M.}
\newblock Ore's conjecture for k=4 and {G}r{\"o}tzsch {T}heorem.
\newblock Manuscript.

\bibitem{hyper}
{\sc Postle, L., and Thomas, R.}
\newblock Hyperbolic families and coloring graphs on surfaces.
\newblock {\em ArXiv e-prints 1609.06749\/} (Sept. 2016).

\bibitem{tw-klein}
{\sc Thomas, R., and Walls, B.}
\newblock Three-coloring {K}lein bottle graphs of girth five.
\newblock {\em J. Combin. Theory Ser. B 92\/} (2004), 115--135.

\bibitem{thom-torus}
{\sc Thomassen, C.}
\newblock Gr\"otzsch's 3-color theorem and its counterparts for the torus and
  the projective plane.
\newblock {\em J. Combin. Theory Ser. B 62\/} (1994), 268--279.

\bibitem{Tho3list}
{\sc Thomassen, C.}
\newblock 3-list coloring planar graphs of girth 5.
\newblock {\em J. Combin. Theory Ser. B 64\/} (1995), 101--107.

\bibitem{thom-surf}
{\sc Thomassen, C.}
\newblock The chromatic number of a graph of girth 5 on a fixed surface.
\newblock {\em J. Combin. Theory Ser. B 87\/} (2003), 38--71.

\bibitem{ThoShortlist}
{\sc Thomassen, C.}
\newblock A short list color proof of {G}rotzsch's theorem.
\newblock {\em J. Combin. Theory Ser. B 88\/} (2003), 189--192.

\bibitem{walls-enum}
{\sc Walls, B.}
\newblock {\em Coloring girth restricted graphs on surfaces}.
\newblock PhD thesis, Georgia Institute of Technology, 1999.

\bibitem{Youngs}
{\sc Youngs, D.}
\newblock 4-chromatic projective graphs.
\newblock {\em Journal of Graph Theory 21\/} (1996), 219--227.

\end{thebibliography}

\end{document}